\documentclass[11pt]{amsart}
\usepackage{amsfonts}

\newtheorem{theorem}{Theorem}[section]
\newtheorem{lemma}[theorem]{Lemma}

\newtheorem{corollary}[theorem]{Corollary}
\newtheorem{proposition}{Proposition}[section]

\theoremstyle{definition}
\newtheorem{problem}[theorem]{Problem}
\newtheorem{remark}[theorem]{Remark}
\newtheorem{definition}[theorem]{Definition}

\newcommand{\MC}{\mathcal{MC}}
\newcommand{\NP}{N}
\newcommand{\NN}{P}

\newcommand{\mhat}{{\mathfrak{m}}}
\newcommand{\ahat}{{\mathfrak{a}}}

\newcommand{\Khat}{U_+}
\newcommand{\Utilde}{\widetilde{U}}
\newcommand{\Ktilde}{\widetilde{K}}
\newcommand{\Kcheck}{\Utilde_+} 

\newcommand{\Ubar}{U}
\newcommand{\UP}{U^\prime}

\newcommand{\Kbar}{K}

\newcommand{\utilde}{\tilde{\mathfrak{u}}}

\newcommand{\kcheck}{\utilde_+}

\newcommand{\shat}{\sigma}
\newcommand{\ub}{\mathfrak{u}}
\newcommand{\ubar}{\mathfrak{u}}
\newcommand{\kbar}{\mathfrak{k}}
\newcommand{\pbar}{\mathfrak{p}}

\newcommand{\kb}{\mathfrak{k}}

\newcommand{\up}{{\mathfrak{u}^\prime}}
\newcommand{\kp}{{\mathfrak{k}^\prime}}
\newcommand{\pp}{{\mathfrak{p}^\prime}}

\newcommand{\ptilde}{\tilde{\mathfrak{p}}}

\newcommand{\tbar}{\tau}

\newcommand{\KP}{K^\prime}
\newcommand{\Ad}{\textup{Ad}}
\newcommand{\khat}{\mathfrak{u}_+}
\newcommand{\phat}{\mathfrak{u}_-}

\newcommand{\azero}{\alpha_0^{++}}
\newcommand{\aminus}{\alpha_1^{+-}}
\newcommand{\aplus}{\alpha_1^{--}}

\newcommand{\real}{{\mathbb R}}
\newcommand{\rt}{\tilde{\rho}}
\newcommand{\VV}{\mathfrak{G}}
\newcommand{\HH}{{\mathbb H}}
\newcommand{\OO}{{\mathbb O}}
\newcommand{\imo}{\textup{Im} {\mathbb O}}
\newcommand{\hh}{\mathcal{H}} 
 
\newcommand{\DD}{{\mathbb D}_+}
\newcommand{\DC}{{\mathbb D}_-}
\newcommand{\cc}{{\mathbb C}}

\newcommand{\bbar}{\left[ \begin{array}}
\newcommand{\ebar}{\end{array} \right] }
\newcommand{\bdm}{\begin{displaymath}}
\newcommand{\edm}{\end{displaymath}}
\newcommand{\beq}{\begin{equation}}
\newcommand{\beqa}{\begin{eqnarray}}
\newcommand{\beqas}{\begin{eqnarray*}}
\newcommand{\eeq}{\end{equation}}
\newcommand{\eeqa}{\end{eqnarray}}
\newcommand{\eeqas}{\end{eqnarray*}}
\newcommand{\dd}{\textup{d}}

\begin{document}

\title[Grassmann geometries and reflective submanifolds]{Grassmann geometries in infinite dimensional homogeneous spaces and an application to  reflective submanifolds}

\author{David Brander}
\address{Department of Mathematics,\\ Faculty of Science, Kobe University,\\1-1 Rokkodai, Nada-ku,\\ Kobe 657-8501, Japan}

\email{brander@math.kobe-u.ac.jp}

\begin{abstract}
Let $U$ be a real form of a complex semisimple Lie group, 
and $(\tau,\, \sigma)$ a pair of commuting involutions on $U$.
This data  corresponds to a reflective submanifold
of a symmetric space, $U/K$.   We 
define an associated integrable system, and
describe how to produce solutions from curved flats.
This gives many new examples of  submanifolds as integrable
systems.

The solutions are shown to correspond to various special submanifolds,
 depending on which homogeneous space, $U/L$, one projects
to. We apply the construction to a question which generalizes, to the context of
reflective submanifolds of arbitrary symmetric spaces,
 the problem of isometric immersions of
space forms with negative extrinsic curvature and flat normal bundle.
For this problem, we prove that the only cases where local solutions 
exist are the previously
known cases of space forms, in addition to our new example of
constant curvature Lagrangian immersions
into complex projective and complex hyperbolic spaces. We also prove non-existence
of global solutions in the compact case. 

For other reflective submanifolds, lower dimensional solutions exist,
and can be described in terms of Grassmann geometries. We consider one example
in detail, associated to the group $G_2$, obtaining a special 
class of surfaces in $S^6$.

\end{abstract}

\keywords{Integrable systems; Loop groups; Symmetric spaces; Reflective submanifolds} 
\subjclass[2000]{Primary 37K10, 37K25, 53C42, 53B25; Secondary 53C35}

\maketitle

\section{Introduction}
The primary aim of this paper is to understand an integrable system characterized 
by a loop group associated to any pair of commuting
involutions on a real semisimple Lie group. This generalizes a system which
was introduced by Ferus and Pedit \cite{feruspedit1996} to study isometric 
immersions of space forms. From this generalization,  new examples of 
special submanifolds as integrable systems can be obtained. 
The application which we will focus on is a natural
analogue, pertaining to an arbitrary reflective submanifold in a symmetric space,
of the problem of isometric immersions of space forms with negative
extrinsic curvature and flat normal bundle.

\subsection{Background}
A \emph{reflective submanifold} $N$,
 of a connected Riemannian manifold $\NN$,
is a totally geodesic symmetric submanifold. Reflective submanifolds
of Riemannian symmetric spaces, $\NN = \Ubar/\Kbar$,
 were classified by DSP Leung \cite{leung1974, leung1975,leung1979}. 
Each can be characterized by a pairwise symmetric Lie algebra,
 $(\ubar, \tbar, \shat)$, where $\tbar$  and $\shat$ are commuting
involutions, and $\tbar$ determines the symmetric space $\Ubar/\Kbar$.
 Let $\ubar = \kbar \oplus \pbar = \khat \oplus \phat$ 
be the canonical decompositions of the Lie algebra, $\ubar$, induced by $\tbar$ and $\shat$ 
respectively.
Then $N$ is the totally geodesic submanifold through the origin, $\Kbar$, 
obtained by exponentiating 
$\pp := \pbar \cap \phat$.  Clearly there is a second (semi-Riemannian) symmetric space 
involved here, $\Ubar/\Khat$, associated to $\shat$.

As part of the project of classifying \emph{symmetric} submanifolds
of symmetric spaces (see \cite{naitoh1986, bent} for this story), H Naitoh
was led to the use of \emph{Grassmann geometries}, a concept 
defined by Harvey and Lawson \cite{harveylawson}.

If $\VV$ is an arbitrary subset of the Grassmann bundle of tangential $s$-planes
over $\NN$, that is, $Gr_s(T\NN) = \cup_{x \in \NN} Gr_s(T_x\NN)$,  then a
\emph{$\VV$-submanifold}, $N$, of $\NN$ is an $s$-dimensional connected 
submanifold such that $T_xN$ is an element of  $\VV$ for each $x \in N$.
 The set of 
such submanifolds is called the \emph{$\VV$-geometry}.  
For us, $\NN$ will be a homogeneous space, $G/H$,
and $\VV$ an orbit of the
action of the identity component, $I_0(G)$, of $G$ on  $Gr_s(T\NN)$.

A $\VV$-geometry
is called \emph{strongly curvature invariant} if, for each $s$-plane $W \in \VV$,
both $W$ and its orthogonal complement $W^\perp$ are invariant under the 
curvature tensor of $\NN$.  
For each strongly curvature invariant $\VV$-geometry of a symmetric
space, $\Ubar/\Kbar$,
there is precisely one reflective submanifold, up to isometries of 
$\Ubar/\Kbar$, and vice versa.  We will be interested in certain
 \emph{non}-totally geodesic $\VV$-submanifolds. For this,
 a result of Naitoh,  
 \cite{naitoh1993, naitoh1998, naitoh2000I, naitoh2000II},
is helpful (see \cite{bent}): 
 the following list contains all strongly curvature invariant $\VV$-geometries
 which admit non-totally geodesic $\VV$-submanifolds,
for the case that $\NN$ is a simply connected,
irreducible symmetric space of compact type:
\begin{enumerate}
\item[(i)]
$k$-dimensional submanifolds of the sphere $S^n$, $0<k<n$;
\item[(ii)]
$k$-dimensional complex submanifolds of the complex projective
space $\cc P^n$, $0<k<n$;
\item[(iii)]
$n$-dimensional totally real submanifolds of $\cc P^n$, $n \geq 2$;
\item[(iv)]
 $2n$-real-dimensional totally complex submanifolds of the quaternionic
projective space $\HH P^n$, $n \geq 2$;
\item[(v)]
The geometries associated with irreducible symmetric $R$-spaces;
 \end{enumerate}
There is  a precisely analogous list for the non-compact case.

\subsection{Problem Statement}
If $E_1$ and $E_2$ are vector bundles over a manifold $M$, with connections
$\Gamma_1$ and $\Gamma_2$ respectively, then we will call the bundle connection
pairs $(E_1,\Gamma_1)$ and  $(E_2,\Gamma_2)$ isomorphic if there exists a vector
bundle  isomorphism $\phi: E_1 \to E_2$ such that
 $\phi^* \Gamma_2 = \Gamma_1$. 
Suppose $M$ is an immersed Riemannian submanifold of $N_1$, and
$f: M \to N_2$ is an isometric  immersion, where 
$\textup{Dim}(N_1) = \textup{Dim}(N_2)$. We will say that $f$ is
 \emph{normal curvature preserving} if the associated normal bundle connection
pairs are isomorphic.

For a Riemannian manifold $M$, let $M_R$ denote the same manifold with the
metric scaled by a factor of $R>0$.  If $\NN := \Ubar/\Kbar$, is a symmetric 
space, and $N := \exp (\pp)$ is a totally geodesic submanifold 
through the origin, $0 := \Kbar$,  whose tangent space
at $0$ is $\pp$, denote the associated orbit on 
$Gr_{\textup{Dim}(\pp)}(T \, \NN)$ by $\VV_\pp$.
%--------------------
\begin{problem}  \label{problem}
Suppose given a simply connected, immersed, reflective submanifold
 $N := \exp(\pp)$ of a semisimple
Riemannian symmetric space $\NN$. Thus, $N_R$ is a reflective submanifold
of $\NN_R$.
Does there exist a (local or global) isometric immersion of $N_R$ into $\NN$ as a 
normal curvature preserving $\VV_\pp$-submanifold?
More specifically, we consider this problem for 
\begin{enumerate}
\item[(i)]
$R>1$, if $\NN$ is of compact type;
\item[(ii)]
$R<1$, if $\NN$ is of non-compact type. 
\end{enumerate}
\end{problem}
In the case that $\NN$ has constant curvature,  these ranges for $R$ are equivalent to the 
requirement that $N_R$ has negative extrinsic curvature.  On the other hand,
for \emph{positive} extrinsic curvature, one has umbilic hypersurfaces $N_R$, so
there is no issue there.

Problem \ref{problem}  generalizes the questions
of whether there is a local isometric embedding with flat normal bundle of a sphere $S_R^k$, of radius $\sqrt{R}$ into $S^n$, for $R>1$, and the analogue for a hyperbolic space $H_R^k$ into $H^n$ for $R<1$.  It is known for these special cases that solutions exist if and only if $n \geq 2k-1$, and, moreover, that the integrability condition is an 
\emph{integrable} system of nonlinear PDE, solutions of which can be obtained by solving a collection of commuting linear differential equations \cite{akns1,
terngtenenblat, terng1980, feruspedit1996}.
For the sphere, even though local solutions exist, there is no complete
global solution \cite{moore1}. 
For the hyperbolic case, $H^k_R \to H^{n}$, for $R<1$, global non-immersibility is also 
conjectured to be the case (usually stated in codimension $k-1$). It
is equivalent to the complete isometric immersion problem with flat normal 
bundle of hyperbolic space into Euclidean space.  
Non-immersibility for $H^2$  into $E^3$ was proven by Hilbert
\cite{hilbert1}. However, in the simply connected case, the problem
 remains open for higher dimensions.

\subsection{Approach}

We consider a generalization of a Grassmann submanifold: given 
a subset, $\VV$, of  $Gr_s(T\NN)$, a \emph{sub-$\VV$-submanifold}, or \emph{$\VV$-compatible}
submanifold is an $r$-dimensional submanifold of $\NN$, $r \leq s$, each of whose tangent
spaces is \emph{contained} in an element of $\VV$. 
The basic idea is,
given two homogeneous spaces, $G/H$ and $G/K$, where $K$ is a subgroup of $H$, 
we construct \emph{sub}-$\VV$-submanifolds in $G/K$ which project to
$\VV^\prime$-submanifolds of $G/H$, for certain geometries $\VV$ and $\VV^\prime$.
In  Section \ref{concremarks} below, we explain that, in this article,
 we are doing precisely that,
in an infinite dimensional homogeneous space.

We study $\VV_{-1}^1$-compatible submanifolds, 
where $\VV_{-1}^1$ is a particular Grassmann geometry
of finite dimensional
submanifolds of an infinite dimensional homogeneous space $\hh/\hh^0$.
Here $\hh$ is a subgroup of the loop group $\Lambda G$
of maps from the unit circle into a complexification, $G$, of $\Ubar$.
The subgroup $\hh$ is defined as the fixed point set
of three involutions of $\Lambda G$, associated to the pairwise symmetric
Lie algebra $(\ubar, \tbar, \shat)$. 
 This is a modification and generalization 
of the set-up used in \cite{feruspedit1996}.

In Section \ref{lc} we outline how to produce $\VV$-compatible submanifolds
for a more general class of Grassmann geometries, which contains ours,
using a generalized Krichever-Dorfmeister-Pedit-Wu (KDPW) method.
This reduces the problem to one of producing simpler objects, which
are then susceptible to the  standard Adler-Kostant-Symes (AKS) theory,
 \cite{adlervanmoerbeke1, adlervanmoerbeke2, kostant1979, symes1980, burstallpedit}.
 These methods  are well-known to experts in the field.

In Section \ref{construction} we describe the three-involution loop group
and the $\VV_{-1}^1$-compatible submanifolds, and work out the basic equations 
satisfied by them. In general, a $\VV_{-1}^1$-compatible immersion 
takes values in the subgroup consisting
of loops which have holomorphic extensions to the punctured plane $\cc^*$.
In our case, for a fixed value of the loop parameter $\lambda \in \real^*$, the
solutions are special (immersed) submanifolds of the homogeneous space
$\frac{\Ubar}{\Kbar \cap \Khat}$, satisfying certain conditions.

%--------------------------
\subsection{Results}
In Section \ref{interpretationsection} we study the projections of these
immersions to $\Ubar/\Kbar$, and then to $\Ubar/\Khat$.  
For the first case, we show, Theorem \ref{firstprojectionthm},
that, if the dimension of the submanifold is equal to the dimension of $\pp$,
 then regular solutions
are precisely the solutions to Problem \ref{problem}, for the case $R>1$ (irrespective
of whether the symmetric space $\NN$ is compact).

For the projections to $\Ubar/\Khat$, we show, Proposition \ref{flatprop},
 that these are
curved flats, that is maps which are, at each point, tangent to a flat of 
$\Ubar/\Khat$. We can then deduce that there is no complete immersion 
(other than a curve) in the compact case:
%------------------------------
\begin{corollary} \label{flatpropcor2}
Let $N$ and $\NN$ be as in Problem \ref{problem}, for the case $\NN$ \emph{compact} and $R>1$, 
$\Sigma_R  \subset N_R$ an open submanifold, 
and $f:\Sigma_R \to \NN$  an isometric normal curvature preserving $\VV_\pp$-immersion.
If the symmetric space $N$ is not flat, then $\Sigma_R$ cannot be the whole of $N_R$.
\end{corollary}

%-----------------
In Section \ref{regularitysection} we study the regularity of the solutions
obtained via the KDPW method, with the following  result, which is basic
to  all that follows:
%--------------------------------------
\begin{theorem} \label{theorem1}
Let $N := \exp(\pp) \subset \NN =\Ubar/\Kbar$ be a reflective submanifold,
where the associated second 
symmetric space $\Ubar/\Khat$ is Riemannian, and let $R>1$. 
Then there exists a local isometric immersion, $\NP_R \to \NN$, as
a normal curvature preserving $\VV_\pp$-submanifold, if and only if 
$\textup{dim} (\NP) \leq \textup{Rank} (\Ubar/\Khat)$.
\end{theorem}
Note that the assumptions that $R>1$ and $\Ubar/\Khat$ is Riemannian
 are essential for the ``only if''
direction.

\subsubsection{The Compact Case}
In Section \ref{compactsection} we check the rank condition for all relevant cases
and conclude:
%------------------------------------
\begin{theorem} \label{compactthm}
The following list contains the geometric interpretations of all possible local solutions
 to Problem \ref{problem}
for the case $R>1$ and  $\NN$ is a simply connected, compact, irreducible, 
Riemannian symmetric
space.  In all cases, local solutions exist and can be constructed by
the AKS-KDPW scheme. In all cases where $\textup{Dim}(N_R) >1$,
 there is no solution which is geodesically
complete.
\begin{enumerate}
\item[(i)]
$N_R=S_R^k$ is an isometric immersion with flat normal bundle from a $k$-sphere of radius
$\sqrt{R}$ into the unit sphere $S^n$, with $0< k \leq (n+1)/2$, and $n \geq 2$.
\item[(ii)]
$N_R=S_R^n$ is an  isometric totally real immersion from an $n$-sphere of radius $\sqrt{R}$
into complex projective space $\cc P^n$, with $n \geq 2$.
\end{enumerate}
\end{theorem}

\subsubsection{The Non-Compact Case}
In Section \ref{noncompactsection} we consider Problem \ref{problem} for the
non-compact case. Here the second symmetric space, $\Ubar/\Khat$ is non-Riemannian,
which makes it less straightforward to obtain non-existence results.  However, we can use
the property of the three involution loop group, explored in \cite{branderdorf, brander2},
that a solution for one problem for a value of the loop parameter $\lambda$ in $\real^*$ corresponds
to a solution to a different problem for $\lambda$ in $S^1$ or $i \real^*$. This allows us,
in those cases where it is necessary, to equate the problem with one in which 
$\Ubar/\Khat$ is Riemannian.  Since the submanifolds are hyperbolic, we do not,
however, get global non-immersibility, as we did in the compact case.
%--------------------------
 \begin{theorem} \label{noncompactthm}
The following list contains the geometric interpretations of all possible local solutions
 to Problem \ref{problem}
for the case $R<1$ and  $\NN$ is a simply connected, irreducible, 
Riemannian symmetric space of non-compact type. 
 In all cases, local solutions exist and can be constructed by
the AKS-KDPW scheme. 
\begin{enumerate}
\item[(i)]
$N_R=H_R^k$ is an isometric immersion with flat normal bundle from a $k$-dimensional
hyperbolic space of constant curvature $\frac{-1}{R}$
into the standard hyperbolic space  $H^n$, with $0< k \leq (n+1)/2$, and
$n \geq 2$.
\item[(ii)]
$N_R=H_R^n$ is an  isometric totally real immersion from the hyperbolic space 
of  constant curvature $\frac{-1}{R}$ into complex hyperbolic space $\cc H^n$, with
$n \geq 2$. 
\end{enumerate}
\end{theorem}

\begin{remark} It seems an interesting question to ask whether or not there
exist complete immersions in the second case of Theorem \ref{noncompactthm},
given the above-mentioned conjectures concerning the first case.
\end{remark}

\subsection{Other Projections}
In the last section we discuss the case where 
$\textup{Dim}(\pp) > \textup{Rank}(\Ubar/\Khat)$.  In this case, the loop group
still provides special submanifolds, but they are of dimension at most equal to
the rank of $\Ubar/\Khat$, so the projections to $\Ubar/\Kbar$ are deformations of
certain submanifolds of symmetric submanifolds, in fact \emph{sub}-$\VV_\pp$-submanifolds.
There are many of these,
as the list of reflective submanifolds contains over one hundred examples \cite{leung1974}.
We consider one example only, associated to a pair of involutions on $G_2$, and
project the solutions, which are of dimension 2, to the sphere $S^6$. 
The resulting surfaces  are deformations of totally geodesic
complex curves in $S^6$, with the property that the restriction of $TS^6$
to the surface decomposes canonically into three 2-dimensional sub-bundles,
each of which is invariant under the almost complex structure of $S^6$.

\subsection{Concluding Remarks} \label{concremarks}
The $\VV_{-1}^1$-compatible submanifolds  here are defined as maps
into $\hh/\hh^0$, where $\hh^0 = \Kbar \cap \Khat$. As such, they are only
\emph{sub}-$\VV_{-1}^1$-submanifolds, and the $\VV_{-1}^1$-geometry on this space
consists of integral submanifolds for a distribution.  However, the final
object we are interested in is the associated family of immersions into
$U/\Kbar$.  This can be thought of as a true $\VV_{-1}^1$-submanifold
in $\hh/\Kbar$, where this geometry is defined analogously to that on $\hh/\hh^0$,
and it is easy to verify that,
 on $\hh/\Kbar$, the $\VV_{-1}^1$-geometry is not a distribution. 

The $\VV_{-1}^1$-compatible submanifolds associated to the three involution 
loop group studied here can, in the real analytic case,
be thought of as families of certain pluriharmonic maps \cite{brander2}.
Moreover, the choices of $R$ complementary to those stated in 
Problem \ref{problem}
can also be obtained from the set-up here. The results for those cases
are different, however, as the symmetric space $\Ubar/\Khat$ is
non-Riemannian.

There are further methods for producing and studying the solutions to the
problems described here, such as dressing, and also studying an associated
linear problem on an algebraic curve. A possible alternative to the KDPW
approach used here would have been the modified AKS construction given in \cite{feruspedit1996}.

For global problems,  the loop group approach is clearly of use in obtaining 
 non-existence results.  On the other hand, the problem of constructing,
or proving existence of, solutions with particular global properties 
is challenging, and more progress in this direction would be interesting.

The basic principal behind this project was to study a naturally 
occurring integrable system
and to interpret the solutions geometrically.  This approach has been pursued 
recently by Terng and collaborators \cite{bdpt2002, terng2002}.

%***********************************************************************
%***********************************************************************

\section{Integrable Systems Constructed from Loop Groups} \label{lc}
In this section we describe a general type of integrable system, and a way
to obtain solutions, using well-established methods. 
We present it in terms of Grassmann geometries, as this fits well
with our applications.
The essential fact we will  need is
the correspondence (\ref{taucor}).

\subsection{Loop Group Grassmann Geometries} 
Let $G$ be a complex Lie group, with 
Lie algebra $\mathfrak{g}$. Let $\Lambda G$ be the group of maps from the 
unit circle $S^1$ into $G$, of a class which includes all loops which extend 
analytically to some neighbourhood of $S^1$, and that makes $\Lambda G$
a Banach Lie group; such as the Wiener topology. $\Lambda G$ includes maps which have
holomorphic extensions to $\cc^*$, and most of the maps we discuss are of this type.

Let $\hh$ be a closed subgroup of $\Lambda G$, and 
 $\hh^0$ be the group of constant loops in $\hh$, namely $\hh \cap G$.
Then the set of left cosets of $\hh^0$, denoted by $\hh/\hh^0$, is a (generally
 infinite dimensional)
homogeneous space on which $\hh$ acts by left multiplication. 
 This induces an action
on the Grassmann bundle $Gr_s \big (T\frac{\hh}{\hh^0} \big )$,
 and so we can talk about associated $\VV$-geometries, where $\VV$
is any orbit of this action.

The Lie algebra $\Lambda \mathfrak{g}$ of $\Lambda G$ consists of the set of 
Fourier series, $\sum_{i= -\infty}^\infty a_i \lambda^i$, $a_i \in \mathfrak{g}$,
 where $\lambda$ is the
$S^1$ parameter, and the convergence condition depends on the topology one chooses.
The Lie algebra, $\textup{Lie}(\hh)$, of $\hh$, is a subalgebra of
$\Lambda \mathfrak{g}$, and the Lie algebra, $\textup{Lie}(\hh^0)$, of $\hh^0$ is
just the set of constant polynomials in $\textup{Lie}(\hh)$.   Hence the tangent space
at the origin, $0 := \hh^0$, of $\hh/\hh^0$ consists of the elements of $\textup{Lie}(\hh)$
whose constant terms are zero, 
\bdm
T_0 \frac{\hh}{\hh^0} = \{ \sum_{i \neq 0} a_i \lambda^i \} \subset \textup{Lie}(\hh).
\edm

Given such a subgroup $\hh$, define, for integers $a<b$,
 $W_a^b$ to be the vector subspace of $T_0 \frac{\hh}{\hh^0}$ given by
\bdm
W_a^b = \{ x \in T_0 \frac{\hh}{\hh^0} ~|~ x = \sum_{i=a}^b a_i \lambda^i \}.
\edm
Now set $\VV_a^b$ to be the distribution given by the 
orbit of $W_a^b$ under the action of $\hh$
on $Gr_{d(b-a)}\big (T \frac{\hh}{\hh^0} \big )$, where $d = \textup{Dim} (\mathfrak{g})$.

The basic object we can construct, using  the techniques described here,
are $\VV_a^b$-compatible
(immersed) submanifolds of $\hh/\hh^0$.
They were studied, essentially from this point of view,  in \cite{branderdorf}, where
they were called connection order $(a,b)$ maps.  This terminology was
grounded in the structure of the Maurer-Cartan form of a frame for 
such an immersion, described below.

Here is a more detailed description:
  let $M$ be a smooth manifold with a fixed base point $x_0$,
 and denote by $\textup{Map}(M,x_0,\Lambda G,I)$ the group of smooth maps
 $M \to \Lambda G$, which map $x_0$ to the identity.   
For a smooth map $f$ from $M$ into a Lie group $\mathcal{G}$,
 we denote by $f^{-1} \dd f$, the pull-back via $f$ of
the Maurer-Cartan form of $\mathcal{G}$.
If $F \in \textup{Map}(M,x_0,\Lambda G,I)$, then the Maurer-Cartan form 
$F^{-1} \dd F$  takes its values in the Lie algebra $\Lambda \mathfrak{g}$,
that is, it has an expansion
\bdm
\alpha = F^{-1} \dd F = \sum _i \alpha_i \lambda ^i,
\edm
where $\alpha_i$ are  $\mathfrak{g}$-valued 1-forms.
For each value of $\lambda$ where $F$ is defined, in particular for all $\lambda \in S^1$,
 we have a smooth map
$F_\lambda :M \to G$, with Maurer-Cartan form $\alpha^\lambda$,
 which must satisfy the integrability condition,
known as the Maurer-Cartan equation,
\beq \label{mce}
\dd \alpha^\lambda + \alpha^\lambda \wedge \alpha^\lambda = 0.
\eeq
Conversely, such a family of 1-forms which satisfies (\ref{mce}) on a simply-connected
manifold $M$ can be integrated to obtain a unique map 
$F \in \textup{Map}(M,x_0,\Lambda G,I)$.

If $\mathcal{H}$ is any subgroup of the loop group $\Lambda G$ then, 
for any  extended integers $a,~b \in {\mathbb Z} \cup \{\pm \infty \}$,
with $a\leq b$,
we define 
\bdm
\MC (M,\hh)_a^b = \{ F \in \textup{Map}(M,x_0, \hh,I) ~| ~ F^{-1} \dd F = \sum_{i=a}^b \alpha_i \lambda^i,~ \alpha_i ~\textup{constant in } \lambda \}.
\edm
The notation ``$\MC$'' is meant to remind us that it is 
the expansion of the Maurer-Cartan form
of $F$, not $F$ itself, which is a polynomial of degree $(a,b)$.

Since postmultiplication by an element which is constant in $\lambda$
has no effect on the connection order, we can analogously define
\bdm
\MC(M,\hh/\hh^0)_a^b.
\edm

If $a$ and $b$ are finite, and 
an element $f \in \MC(M,\hh/\hh^0)_a^b$
is a regular map, $M \to \hh/\hh^0$, 
then $f$ is just a $\VV_a^b$-compatible immersion into $\hh/\hh^0$, mapping
$x_0$ to the origin.
 On the other hand, 
for a fixed value of $\lambda \in \cc^*$, these are maps, $f^\lambda$, from
$M$ into some quotient group  $\frac{U}{\mathcal{H}^0}$, where $U$ is some
subgroup of $G$.  This leads to the interpretations of 
$\VV_a^b$-compatible immersions as families (as $\lambda$ varies)
 of special submanifolds of a homogeneous space.
The fact that the Maurer-Cartan equation (\ref{mce}) holds for all $\lambda$,
imposes some conditions on $f_\lambda$, which define the special submanifold.

\noindent \textbf{Notation:} We sometimes use the same symbol $x$ for an element of
$\hh/\hh^0$ as for a representative in $\hh$.

\subsection{The KDPW Method} \label{dpw}
This is a general method for constructing a $\VV_a^b$-compatible immersion,
where $a<0<b$, from a pair of simpler objects, namely a $\VV_a^0$-compatible
immersion and a $\VV_0^b$-compatible immersion, using a Birkhoff factorization
of the loop group (see \cite{pressleysegal}). This is a factorization for
$x \in \Lambda G$, 
\beq \label{birkhoff}
x = x^+ x^-,
\eeq 
where the loops $x^+$ and $x^-$ have holomorphic extensions to the unit disc
$\DD$ and the complement of its closure $\DC$ respectively. 
Denote the respective subgroups of such elements by $\Lambda^\pm G$,
and, for a subgroup $\hh$, set $\hh^\pm := \hh \cap \Lambda^\pm G$.
The Birkhoff factorization
is valid on an open dense neighbourhood of the identity of $\Lambda G$, and
is unique with the normalization $x^+(0) = I$. Obviously the analogue holds
with $\pm$ interchanged. 

The idea was used by Krichever in \cite{krichever} for solutions of the
sine-Gordon equation. A different variant was used
for harmonic maps into symmetric spaces
by Dorfmeister, Pedit and Wu in \cite{dorfmeisterpeditwu}.
In the case of harmonic maps, it reduced the problem to specifying certain
holomorphic data, and has since been used successfully to produce numerous 
examples of constant mean curvature surfaces. This is known as the DPW
method. Dorfmeister and the author investigated the
general applicability of the basic idea in  \cite{branderdorf}.

If $x$ is an element of a subgroup $\hh$ of $\Lambda G$, and 
 the factors $x^\pm$ in (\ref{birkhoff}) are also in $\hh$, then
we say $\hh$ is  Birkhoff decomposable.  
In this case, the generalized KDPW method
described in \cite{branderdorf} gives the following bijection:
\beqa \label{bcor}
F \in \MC(\Omega,\hh/\hh^0)_a^b, ~~ \leftrightarrow  &
 F_+ \in \MC(\Omega,\hh)_1^b, \\
 & F_- \in \MC(\Omega,\hh)_a^{-1},  \nonumber
 \eeqa
which holds at least in a neighbourhood, $\Omega$, of the base point,
 $x_0$, of $M$.
The maps $F_+$ and $F_-$ are simply the unique normalized 
left factors in the left and right 
Birkhoff decompositions $F =F_+G_- = F_-G_+$. 

The version which we will use, is the case that $\hh$ is Birkhoff
decomposable and we consider a further subgroup, $\hh_\tau$,
where $\tau$ is an involution of the second kind of $\hh$, meaning it takes
$\hh^\pm \to \hh^\mp$,
and $F \in \MC(M,\frac{\hh_\tau}{\hh^0_\tau})_a^b$:
then it follows that $a = -b$ and,  applying
$\tau$ to both sides of the correspondence (\ref{bcor}), we see that
we must have $F_- = \tau F_+$. One can show that there is
always a bijective correspondence
\beq  \label{taucor}
 F \in \MC(\Omega,\frac{\hh_\tau}{\hh^0_\tau})_{-b}^b
 ~~ \leftrightarrow  ~~ F_+ \in \MC(\Omega,\hh)_1^b, 
  \eeq
for some neighbourhood $\Omega$ of $x_0$.
If $\hh_\tau$  is a subgroup of $\Lambda U$, where $U$ is a compact real
form of $G$, then 
 we also have a global statement in the right to left direction, because
$F$ on the left hand side of (\ref{taucor}) is constructed pointwise from
 an Iwasawa factorization.

%***********************************************************

\subsection{The AKS Theory and Finite Type Solutions}
One useful consequence of the correspondence (\ref{bcor}) is
that the standard AKS theory can be used to generate many solutions
which are computed by solving a collection of ODE's on a 
finite dimensional vector space, so-called finite type solutions.
In general, this method applies to $\VV_a^{-1}$-compatible
and $\VV_1^b$-compatible immersions, and one can then apply the
correspondence $(\ref{bcor})$ to these. 
A short description can be found in  \cite{burstallpedit}.
We will not use this construction here, but wish to point out that this
is one means to generate many solutions for the problems studied below.

%**************************************************************
\section{The Three Involution Loop Group} \label{construction}
\subsection{Definition of the Group} \label{groupdefinition}
Let $G$ be a complex semisimple Lie group, $\tbar$, $\shat$ a pair
of commuting involutions of $G$, and $\rho$ a conjugation of $G$
which commutes with $\tbar$ and $\shat$.   Let $\Ubar := G_\rho$ be the
real form determined by $\rho$, with corresponding Lie algebra $\ubar$.
Extend the involutions to $\Lambda G$ by the rules:
\beqa \label{symrho}
(\rho X)(\lambda) = \rho(X(\bar{\lambda})),\\
(\shat X)(\lambda) = \shat(X(-\lambda)), \label{symshat}\\
(\tbar X)(\lambda) = \tbar(X(-1/\lambda)), \label{symtau}
\eeqa
and consider the subgroup fixed by all three involutions:
\bdm
\hh = \Lambda G_{\rho \tbar \shat}.
\edm
It is easy to show that 
if $\Lambda G$ is any Birkhoff decomposable group, for example
if $G$ is complex semisimple, then so is the subgroup defined as the fixed
point set of any finite
number of commuting finite order automorphisms of the first kind, which include $\rho$ and
$\shat$ defined here. Since $\tau$ is of the second kind,
 we are in the realm of the second version of the KDPW method and, in particular,
can use the correspondence (\ref{taucor}).

Now consider a $\VV_{-1}^1$-compatible immersion in $\hh/\hh^0$, that is a regular map   
\bdm
 F \in \MC(M,\hh/\hh^0)_{-1}^1.
\edm 

For $\lambda \in \real^*$, $F$ is a map from $M$ into
$\Ubar/\Ubar_{\tbar \shat}$,  since
$\Lambda G_{\rho \tbar \shat}^0 = \Ubar_{\tbar \shat}$. We can also project to obtain 
maps into the symmetric spaces $\Ubar/\Ubar_{\tbar}$ and $\Ubar/\Ubar_{\shat}$, or more
generally, into any homogeneous space $\Ubar/H$, where $\Ubar_\tbar \cap \Ubar_{\shat} \subset H$.

In principle, we can construct $F$  from  elements
\bdm
F_+ \in \Lambda G_{\rho \shat}(M)_1^1,
\edm
by the correspondence (\ref{taucor}). 
The  maps $F_+$ are families of curved flats  into $\Ubar/\Ubar_{\shat}$.

%----------------------------------------------------
\subsection{Curved Flats} \label{curvedflatssection}
Here we summarize the definitions and relevant results from \cite{feruspedit1996II},
where Ferus and Pedit defined and studied curved flats.
\begin{definition}
If $M$ is a manifold and $N=U/K$ is a (semi-Riemannian) symmetric space, then
$f: M \to N$ is called a \emph{curved flat} if the pull-back via $f$ of the curvature
form on $U/K$ is zero. A curved flat is \emph{regular} if it is immersive everywhere.
\end{definition}
 If $\mathfrak{u} = \mathfrak{k} \oplus \mathfrak{p}$ is the canonical decomposition,
and the tangent bundle to $N$ is identified with $\mathfrak{p}$ via the differential of
left translation from the origin, then the tangent space at each point of a curved flat is
an Abelian subspace of $\mathfrak{p}$.  Hence if the symmetric space is Riemannian, the maximum
possible dimension for a regular curved flat is $\textup{Rank}(U/K)$.

If $M$ is simply connected, then curved flats can be characterized as elements of 
\bdm
\MC(M, \widehat{\hh}/\widehat{\hh}^0)_0^1,
\hspace{1cm} \widehat{\hh} := \Lambda G_{\rho \sigma},
\edm
where $\rho$ is the conjugation determining the real form $U$ of a
complexification, $G$, of $U$, and $\rho$ is extended to 
$\Lambda G$ by $(\rho X)(\lambda) = \rho (X(\bar{\lambda})$,
  $\sigma$ is extended by the formula
$(\sigma X)(\lambda) = \sigma(X(-\lambda))$,
and $K$ is the fixed point set of $\sigma \big |_U$. 
To see this, let $F: M \to U$ be a frame for a curved flat $f$. Note that a global frame
may not exist, but one can show that the map into the homogeneous space 
$\widehat{\hh}/\widehat{\hh}^0$ we are about
to construct is nevertheless defined.  The Maurer-Cartan form of $F$ has the 
$\mathfrak{k} \oplus \mathfrak{p}$ decomposition
\bdm
\alpha = F^{-1}\dd F = \alpha_{\mathfrak{k}} + \alpha_{\mathfrak{p}}.
\edm
The curved flat condition is given by
\bdm
\alpha_{\mathfrak{p}} \wedge \alpha_{\mathfrak{p}} = 0.
\edm
It is simple to check that this condition, together with the integrability condition
for $\alpha$,  namely the Maurer-Cartan equation (\ref{mce}), are equivalent to
the assumption that the family of 1-forms,
\bdm
\alpha^\lambda := \alpha_{\mathfrak{k}} + \lambda \alpha_{\mathfrak{p}},
\edm
satisfies the Maurer-Cartan equation for all values of $\lambda$.
For $\lambda \in \real^*$, we can integrate $\alpha^\lambda$ on $M$ to obtain a frame
for a curved flat.  The associated family is an element of 
$\MC(M, \widehat{\hh}/\widehat{\hh}^0)_0^1$.
Note that, if $M$ is contractible, so that a global frame exists,
 we can gauge away the $\alpha_{\mathfrak{k}}$ term, to get a unique (up to the normalization
 point $x_0$) element of
 $\MC(M, \widehat{\hh})_1^1$.

The results which are of most  interest to us are:
\begin{theorem} \cite{feruspedit1996II} \label{curvedflatthm}
\begin{enumerate}
\item[(i)]
If $\textup{Rank}(U/K) = r$, then the AKS scheme can be used to locally construct infinitely
many $r$-dimensional curved flats into $U/K$, which are, at each point,
 tangent to a Cartan subalgebra of $\mathfrak{p}$.
 \item[(ii)]
If $U/K$ is Riemannian, then a curved flat in $U/K$ is intrinsically flat.
 \end{enumerate}
\end{theorem}

The first statement implies that, on $V$,  the $\VV_{-1}^1$-compatible maps, 
$F \in \MC(M,\hh/\hh^0)_{-1}^1$,
constructed via the KDPW scheme described in Section {\ref{dpw},
will also be immersions into $\Ubar/(\Ubar_{\tbar} \cap \Ubar_{\shat})$, for $\lambda \in \real^*$.
Regular projections to $\Ubar/\Ubar_{\tbar}$ and $\Ubar/\Ubar_{\shat}$ will be discussed
below.

\subsection{The Basic Equations Satisfied by $\VV_{-1}^1$-Compatible Immersions}
\label{equationssection}
Let us return to the loop group defined in Section \ref{groupdefinition}.
Let $F$ be a local frame for  an element of  $\MC(M,\hh/\hh^0) _{-1}^1$,
 that is, $F \in \MC(\Sigma,\hh)_{-1}^1$, for some contractible open set
 $\Sigma$ containing the base point $x_0$.
$F$ takes its values in $\Ubar := G_\rho$ for $\lambda \in \real^*$, so, to simplify notation,
we will always assume $\lambda \in \real^*$ in the following discussion.
As in the introduction, let 
\bdm
\ubar = \kbar \oplus \pbar = \khat \oplus \phat
\edm
be the canonical decompositions associated to $\tbar$ and $\shat$ respectively.
We have the orthogonal (with respect to the Killing-form of $\mathfrak{g}$) decomposition
\bdm
\ub = \ub^{++} \oplus \ub^{+-} \oplus \ub^{-+} \oplus \ub^{--}, 
\edm
where
\beqas
\ub^{++} =   \kb \cap \khat =: \kp, && \hspace{1cm}
 \ub^{+-} =  \kb \cap \phat,\\
 \ub^{--} =  \pbar \cap \phat =: \pp, && \hspace{1cm}
 \ub^{-+} = \pbar \cap \khat.
 \eeqas
Let $\Omega(W)$ denote the set of real-valued 1-forms on a manifold $W$.
%***************************************************************
\begin{lemma} \label{conditionlemma}
\begin{enumerate}
\item[(i)]
If $F$ is any element of $\MC(\Sigma,\hh)_{-1}^1$, then the Maurer-Cartan form 
$\alpha^\lambda= F^{-1}\dd F$ has the expansion
\beq  \label{alphamcf}
\alpha^\lambda= \alpha_0^{++} + \alpha_1^{+-} (\lambda - \lambda^{-1}) +  \alpha_1^{--}(\lambda + \lambda^{-1}),
\eeq
where
\bdm
\nonumber \alpha_0^{++} \in \ub^{++} \otimes \Omega(M), \hspace{.5cm}
 \alpha_1^{+-} \in   \ub^{+-} \otimes \Omega(M),   \hspace{.5cm}
  \alpha_1^{--} \in \ub^{--}\otimes \Omega(M).
\edm  
For any fixed value of $\lambda \in \real^*$, we have, in addition to the integrability of $\alpha$,
the equation
\beq \label{lemcond}
\dd \azero + \azero \wedge \azero = -4\aplus \wedge \aplus.
\eeq

%-------------------
\item [(ii)]
Conversely, suppose  $\alpha^\lambda$ is a family of 
$\ubar$-valued 1-forms on $\Sigma$, of the form (\ref{alphamcf}), for $\lambda \in \real^*$,
 which is integrable for all $\lambda$. Then there exists a unique map 
 $F \in \MC(\Sigma,\hh)_{-1}^1$ whose Maurer-Cartan form is $\alpha^\lambda$.
\end{enumerate}
\end{lemma}
%**********************
\begin{proof}
\begin{enumerate}
\item[(i)]
By definition, $\alpha$  has the expansion
\beqas
\alpha^\lambda= \sum_{i=-1}^1 \, (\alpha_i^{++} + \alpha_i^{+-} + \alpha_i^{-+} + \alpha_i^{--})\lambda ^i ,\\
\alpha_i^{++} \in \ub^{++} \otimes \Omega(M), \hspace{.5cm}
 \alpha_i^{+-} \in   \ub^{+-}\otimes \Omega(M),\\
  \alpha_i^{-+} \in \ub^{-+}\otimes \Omega(M), \hspace{.5cm}
   \alpha_i^{--} \in \ub^{--}\otimes \Omega(M).
\eeqas   
Now 
\bdm
\tbar \alpha^\lambda=   \sum_{i=-1}^1 \, (\alpha_i^{++} + \alpha_i^{+-} - \alpha_i^{-+} - \alpha_i^{--})(-\lambda )^{-i},
\edm
and $\alpha$ is fixed by $\tbar$, from which it follows that 
\bdm
\alpha^\lambda= \alpha_0^{++} + \alpha_0^{+-} + (\alpha_1^{++} + \alpha_1^{+-})(\lambda - \lambda^{-1})
  +  (\alpha_1^{-+} + \alpha_1^{--})(\lambda + \lambda^{-1}).
  \edm
Then 
\bdm 
\shat \alpha^\lambda=  \alpha_0^{++} - \alpha_0^{+-} -(\alpha_1^{++} - \alpha_1^{+-})(\lambda - \lambda^{-1})
  -  (\alpha_1^{-+} - \alpha_1^{--})(\lambda + \lambda^{-1}),
\edm
from which we conclude that $\alpha$ has the form given by (\ref{alphamcf}).

Now let us expand the integrability condition
\beq \label{amce}
\dd \alpha^\lambda+ \alpha^\lambda\wedge \alpha^\lambda= 0.
\eeq
This equation must hold for all values of $\lambda$, and it is straightforward to verify
that this is equivalent to the following five equations
\beqa  \label{firstcond}
\dd \azero + \azero \wedge \azero + 2(\aplus \wedge \aplus - \aminus \wedge \aminus) = 0,\\
 \dd \aminus + \azero \wedge \aminus + \aminus \wedge \azero = 0, \label{secondcond} \\
\dd \aplus + \azero \wedge \aplus + \aplus \wedge \azero = 0, \label{thirdcond} \\
 \aminus \wedge \aplus + \aplus \wedge \aminus = 0, \label{fourthcond} \\
\aplus \wedge \aplus = -\aminus \wedge \aminus.  \label{lastcond}
\eeqa
Now, at a fixed value of $\lambda$, the first four of these equations are merely the decomposition of 
(\ref{amce}) into its $\ub^{++}$, $\ub^{+-}$, $\ub^{--}$ and $\ub^{-+}$ components respectively,
and so, pointwise in $\lambda$,  they give no new conditions beyond integrability.  The last, however,
equation (\ref{lastcond}), does imply a further condition, and, given the equation (\ref{firstcond}),
it is equivalent to equation (\ref{lemcond}) of the Lemma.\\

\item[(ii)]
The converse is just the standard fact  that a $\mathfrak{g}$-valued
1-form which satisfies the integrability condition is the Maurer-Cartan form
of a smooth map into a Lie group $G$, unique up to an initial condition. This
is applied for every $\lambda$, and the uniqueness is given by the normalization
at the base point $x_0$ of $M$.
\end{enumerate}
\end{proof}

%**************************************************

\section{Interpretation of the Solutions}  \label{interpretationsection}
Given an element $F \in \MC(M,\hh/\hh^0)_{-1}^1$, for real values of $\lambda$, $F_\lambda$ is a
map $M \to {\Ubar}/(\Kbar\cap \Khat)$, as explained in Section \ref{construction}. 
We will interpret the projections to both $\Ubar/\Kbar$ and $\Ubar/\Khat$, but first,
some more details about reflective submanifolds.

%****************************************** symmetric spaces *******************
\subsection{Strongly Curvature Invariant Immersions Associated to Reflective Submanifolds of Symmetric Spaces}
\label{constructionsection}
In this section we outline some well-known facts. More details can be found in 
\cite{leung1974}.
\subsubsection{Totally Geodesic Submanifolds of Symmetric Spaces}  \label{outline}
Let $\NN = \Ubar/\Kbar$ be a 
(not necessarily Riemannian) connected semisimple symmetric space,
where $\Kbar = \Ubar_{\tbar}$, and
$\ub = \kb \oplus \pbar$ is the canonical decomposition of the Lie algebra of $\Ubar$.
  There is a natural one-to-one correspondence between 
connected totally 
geodesic submanifolds $\NP$ through the origin $0$ of $\NN$ and linear subspaces 
$\pp $ of $\pbar$ which are closed
under the Lie triple product, that is, $[\pp,[\pp,\pp]] \subset \pp$.
 One way to construct
the submanifold $\NP$ from $\pp$ is by setting $\kp:= [\pp,\pp]$, 
$\up := \kp + \pp$,  
defining $\UP$ to be the connected Lie subgroup of $\Ubar$ generated by $\up$, and setting 
$\KP := \UP \cap \Kbar$.  Then $\NP = \UP/\KP$, and
 the tangent space at the origin, $T_0 \NP$, is $\pp$.
The involution $\tau^\prime$ defining $\KP$ is given by $\tau ^\prime:= \tbar|_\up$.
Conversely, given such a totally geodesic submanifold, one can take $\UP$ to be the
largest subgroup of $\Ubar$ which leaves $\NP$ invariant, and set $\KP := \UP \cap K$.

\subsubsection{Reflective Submanifolds} \label{reflective}
Now suppose given a symmetric space $\Ubar/\Kbar$ as above, and suppose given another
involution $\hat{\sigma}$ of $\Ubar$ which commutes with $\tbar$.
Let $\ubar = \khat \oplus \phat$ be the canonical decomposition associated to $\shat$.
Now set 
\bdm
\pp := \phat \cap \pbar,
\edm
 then
$[\pp,[\pp,\pp]] \subset \pp$, and so there is an associated symmetric
 subspace   $(\UP,\KP,\tau^\prime)$ of
$(\Ubar,\Kbar,\tbar)$, with $\up = \kp + \pp$, 
 and a totally geodesic submanifold $\NP = \UP/\KP$ of
$\NN$, as described above.

It will be convenient sometimes to choose a larger symmetric subalgebra than
 the one described in
 Section \ref{outline} to represent $N$;  namely, set 
\bdm
\kp := \kb \cap \khat,  \hspace{1cm} \up := \kp + \pp.
\edm
Clearly $[\kp,\kp] \subset \kp$ and $[\kp,\pp] \subset \pp$, so 
$(\up, \kp, \tau^\prime)$ is a symmetric subalgebra
of $(\ub, \kb, \tbar)$, where $\tau^\prime := \tbar |_\up$,
 and we have a corresponding symmetric subspace
$(\UP,\KP,\tau^\prime)$.
  Since a totally geodesic submanifold through the origin 
of $\Ubar/\Kbar$ is determined
by its tangent space at the origin (in this case $\pp$) it follows that $\UP /\KP$
is the submanifold corresponding to $\pp$ described  in Section \ref{outline}.

Note that in the situation just described we have the orthogonal decomposition
\bdm
\pbar = \pp \oplus \pp^\perp, \hspace{1cm} \pp^\perp =  \pbar \cap \khat,
\edm
and it is easy to check that $\pp^\perp$ is also closed under the Lie triple product.
Thus $\NN$ is a \emph{reflective} submanifold.

If $\Ubar/\Kbar$ is Riemannian then one can also show that, conversely, every
reflective submanifold is associated canonically to a second involution $\shat$
which commutes with $\tbar$, and so may be characterized by the pairwise symmetric
Lie algebra $(\ubar, \tbar, \shat)$.

In the Riemannian case we also have  the orthogonal decompositions
\beqas
\kb = ([\pp,\pp] \cup [\pp^\perp,\pp^\perp]) \oplus [\pp,\pp^\perp],\\
\pbar = \pp \oplus \pp^\perp,
\eeqas
and the relations
\beqas
~[[\pp, \pp],\pp^\perp] = [[\pp,\pp^\perp],\pp] \subset \pp^\perp,\\
~[[\pp^\perp, \pp^\perp],\pp] = [[\pp,\pp^\perp],\pp^\perp] \subset \pp.
\eeqas

\subsubsection{The $\VV_\pp$-Geometry}
Any pairwise symmetric Lie algebra $(\ubar, \tbar,\shat)$
induces a strongly curvature invariant Grassmann geometry, which
we call the \emph{$\VV_{\pp}$-geometry}, defined by the orbit of 
$\pp = \pbar \cap \phat$ under the action of $\Ubar$ on 
$Gr_{\textup{Dim}(\pp)}(T \frac{\Ubar}{\Kbar})$.
A $\VV_{\pp}$-submanifold $N$ of $\NN = \Ubar/\Kbar$ is characterized by the property that
it is $\textup{Dim}(\pp)$-dimensional, and its tangent space at each point is $\pp$, if we identify the tangent bundle of $\NN$ with $\pbar$ via left translation.

%****************************************************************************

\subsection{The Projection to $\Ubar/\Kbar$}  \label{firstprojectionsection}
Let $N= U/K$ be the reflective submanifold associated with the pairwise symmetric
Lie algebra $(\ubar,\tbar,\shat)$, and set $\pp = \pbar \cap \phat$ as before.
Set 
\beq \label{Rlambda}
R_\lambda := \frac{(\lambda + \lambda^{-1})^2}{4}.
\eeq
For simplicity, we consider the spectral parameter in the interval
 $\real_+ := (0,\infty) \in \real$. The following would analogously
 apply to the interval $(-\infty,0)$.
 As $\lambda$ varies through $\real_+$, $R_\lambda$ varies through the
interval $[1,\infty)$, achieving the value 1 only at $\lambda =1$.

%*****************************
\begin{theorem} \label{firstprojectionthm}
\begin{enumerate}
\item[(i)]
Let $F$  be an element of   $\MC(M,\hh/\hh^0)_{-1}^1$, where 
$M$ is simply connected and $\textup{Dim}(M) = \textup{Dim}(\pp)$.
 Let $f_\lambda :M \to \NN$
be the family of maps obtained by projecting  $F_\lambda$, for $\lambda \in \real_+$, and assume that $f_\lambda$ is regular on $M$ for some value of $\lambda$. Then:
%---------------------------
\begin{enumerate}
\item \label{aa}
$f_{\{\lambda = 1\}}(M)$ is an open submanifold of the reflective submanifold 
$N \subset \NN$. 
\item  \label{bb}
Let $M_R$ be the manifold $M$ with the pull-back metric induced from
$f_{\{\lambda = 1\}}: M \to N$, scaled by a factor $R>0$, so that $M_R$ is a
open  submanifold of $N_R$.
Then, for other values of $\lambda \in \real_+$,  
$f_\lambda$ is an
isometric normal curvature preserving $\VV_\pp$-immersion,
 $f_\lambda: M_{R_\lambda} \to \NN$. \\
\end{enumerate}
%---------------------------
\item[(ii)]
Conversely, let $M$ be a simply connected,  immersed, open submanifold of $N$, $R>1$, 
and suppose $f:M_R \to \NN$ is an isometric normal curvature preserving $\VV_\pp$-immersion.
Choose $\lambda _0 \in \real_+$ such that $R= R_{\lambda_0}$. 
Then there exists an element $F \in \MC(M,\hh/\hh^0)_{-1}^1$
such that the projection of $F$ to $\Ubar/\Kbar$ at $\lambda = \lambda_0$ is $f$.
\end{enumerate}
\end{theorem}
%------------------------------
\begin{proof}
\begin{enumerate}
\item[(i)]
The conditionis that need to be satisfied   are essentially
local, so we may as well assume $M$ is contractible and 
take a global frame for $F$, which we also denote by $F$.
The left translation to the identity of the tangent space to the symmetric space $\Ubar/\Kbar$ is just $\pbar$. 
Let $\langle \cdot\, , \, \cdot \rangle$ be the canonical metric on $\pbar$, given by
the Killing form.  Choose a frame for $T_xM$, $e_i$, $i = 1,...,m$,
 which is orthonormal with respect to the pull-back metric. Then
  we have the projection to $\pbar$ of $\alpha$,
\bdm
\pi_{\pbar} F^{-1} \dd F =: \theta ^\lambda = \sum_i \theta^\lambda_i,
\edm 
where $\theta_i^\lambda (e_i) =: \hat{e}_i$ are orthonormal with respect to 
$\langle \cdot\, , \, \cdot \rangle$,
that is, $\theta_i^\lambda$ are the coframe to $e_i$, with indices lowered.
Comparing the definition of $\theta^\lambda$ 
with the form of $\alpha$ given by (\ref{alphamcf}), we can also see that
\bdm
\theta^\lambda =  \aplus(\lambda +\lambda^{-1}),
\edm
which takes its values in $\pp$.
It follows that $f_\lambda$ is a $\VV_\pp$-immersion.

Now the pullback metric is given by
\bdm
\dd s^2_\lambda = \sum_i \langle \theta^\lambda_i,\theta^\lambda_i \rangle.
\edm
 Evidently, for any real value of $\lambda$, the metric
is just a constant multiple of the metric for the immersion obtained at $\lambda = 1$,
that is:
\bdm
\dd s_\lambda^2 = \frac{(\lambda + \lambda^{-1})^2}{4} \dd s_1^2.
\edm
At $\lambda =1$, the Maurer-Cartan form of $F$ is:
\beq  
\alpha^1 = \alpha_0^{++} + \theta^1,
\eeq
which takes its values in $\up = \kp \oplus \pp$.
It follows that, at $\lambda = 1$, $F$ is a frame for a map into the totally geodesic
submanifold $\UP/\KP$, in other words, $f_1$ is an isometric immersion $M \to \UP/\KP$. 
For other values of $\lambda \in \real_+$,
 $f_\lambda$ is an isometric immersion
$M_{R_\lambda} \to \Ubar/\Kbar$.  

Finally, the connection 1-form for $f_\lambda$ is  given by the projection to $\kb$ of the
Maurer-Cartan form of $F_\lambda$, which, from (\ref{alphamcf}) is
\bdm
\pi_{\kb} \alpha^{\lambda} = \alpha_0^{++} + \alpha_1^{+-} (\lambda - \lambda^{-1}),
\edm
where $\alpha_0^{++} \in \kb \cap \khat \otimes \Omega(M)$ and 
 $\alpha_1^{+-} \in \kb \cap \phat \otimes \Omega(M)$.
We have the relations:
\beqa  \label{1strelation}
~[\kb \cap \khat, \pp] \subset \pp, \hspace{1cm}
[\kb \cap \khat, \pp^\perp] \subset \pp^\perp,\\
~[\kb \cap \phat, \pp] \subset \pp^\perp, \hspace{1cm} 
[\kb \cap \phat, \pp^\perp] \subset \pp. \label{2ndrelation}
\eeqa
It follows that the $\kp$ component of the connection, $\alpha_0^{++}$,
 is the sum of the connections on $TM \big |_{M_{R_\lambda}}$ and $T^\perp M \big |_{M_{R_\lambda}}$,
whilst $\alpha_1^{+-} (\lambda - \lambda^{-1})$ is the second fundamental form.

Fixing a choice of bases for $\pp$ and $\pp^\perp$ gives a  bundle isomorphism between
$TM \big |_{M_{R_\lambda}}$ and  $TM \big |_M$ as well as between
$T^\perp M \big |_{M_{R_\lambda}}$ and $T^\perp M \big |_{M}$.
The sum of the tangential and normal connections do not depend on $\lambda$, 
so this 1-form pulls-back identically via this isomorphism.  Since we already
know that the tangent bundle-connection pairs are isomorphic (the metrics 
differ by a constant scalar multiple), it follows that the normal bundle-connection
pairs are also isomorphic.\\

\item [(ii)]
For the converse, we argue locally first, that is, choose a contractible subset
$\Sigma_R$ of $M_R$, which contains the normalization point, so that a frame $F: \Sigma_R \to U$
exists for $f$. 
The Maurer-Cartan form of $F$ has the decomposition
\beqas
F^{-1}\dd F =  \alpha^{++} + \alpha^{+-} + \alpha^{-+} + \alpha^{--} ,\\
\alpha^{\pm \pm } \in \ub^{\pm \pm} \otimes \Omega(M).
 \eeqas   
The  assumption  that $f$ is  a $\VV_\pp$-immersion means that $T^\perp \Sigma_R = \pp^\perp \big |_{\Sigma_R}$, where, as usual, we identify $T\NN$ with $\pbar$. 
Hence the projection
to $\pp^\perp$ of $F^{-1}\dd F$ is zero, that is, $\alpha^{-+} = 0$.
Now define a 1-form 
\beq \label {requiredform}
\alpha^\lambda =  \alpha^{++} + 
    (\lambda -\lambda^{-1}) \frac{\alpha^{+-} }{ \lambda_0 -\lambda_0^{-1}}
    + (\lambda + \lambda^{-1}) \frac{\alpha^{--}} {\lambda_0 + \lambda_0^{-1}} .\\
\eeq
By definition,  $\alpha^{\lambda_0} = F^{-1}\dd F$ and is therefore integrable.
Further, we saw in Lemma \ref{conditionlemma} that $\alpha^{\lambda}$ being
integrable for \emph{all} values of $\lambda$ is equivalent to integrability
at a single value of $\lambda$ plus the equation
\beqa \label{needtoshow5}
\dd \alpha^{++} +  \alpha^{++}\wedge \alpha^{++}  &=& -\frac{4}{(\lambda_0 + \lambda_0)^2}\alpha^{--} \wedge \alpha^{--},\\
&=& -\frac{1}{R}\alpha^{--} \wedge \alpha^{--}.
\eeqa
We also showed  in Lemma \ref{conditionlemma} that 
$\MC(\Sigma_R,\hh/\hh^0)|_{-1}^1$ consists precisely of families, for $\lambda \in \real^*$,
 of  equivalence classes of 
frames $F_\lambda: \Sigma_R \to U$ whose Maurer-Cartan forms have the form
(\ref{requiredform}) and satisfy (\ref{needtoshow5}),
so we only need to show that this is satisfied.

The 1-forms $\alpha^{--}$ and $\alpha^{++}$ have the same interpretation
as in the first part of the proof.
On the other hand, we also assumed that $M$ is an open immersed submanifold of the reflective
submanifold $N$. If $g: M \to N \subset \NN$ is this immersion, and $G: \Sigma \to \Ubar$
a frame, where $\Sigma = \Sigma_R$ as a point set, then the Maurer-Cartan form is
\bdm
G^{-1} \dd G = \beta^{++} + \beta^{--},
\edm
with $\beta^{\pm \pm } \in \ub^{\pm \pm} \otimes \Omega(M)$. There is no $+-$ term,
because $\Sigma$ is a totally geodesic submanifold of $\NN$, and, from the relation
(\ref{2ndrelation}),  the $+-$ term of
the Maurer-Cartan is  the second fundamental form.
Again, $\beta^{--}$ is the coframe for $g$, and $\beta^{++}$ is the sum of the 
connections on $T \Sigma$ and $T^\perp \Sigma$.  

Because $f$ is an isometric normal curvature preserving
immersion, we can choose a bundle isomorphism 
$\phi: T\Sigma_R  \oplus T^\perp \Sigma_R \to
T \Sigma \oplus T^\perp \Sigma$, such that, 
\beq \label{pullbackeqn}
\phi^* \beta^{++} = \alpha^{++}, \hspace{1cm} \sqrt{R} \, \phi^* \beta^{--} = \alpha^{--}.
\eeq
The factor of $\sqrt{R}$ is due to the fact that the metric on $\Sigma_R$ is $R$ times
that on $\Sigma$.
But then the $\kb\cap \khat$ component of the integrability condition for $G^{-1} \dd G$ is
\bdm
\dd \beta^{++} +  \beta^{++}\wedge \beta^{++}  = - \beta^{--} \wedge \beta^{--},
\edm
and this, combined with (\ref{pullbackeqn}) is equivalent to (\ref{needtoshow5}).
Hence the 1-form $\alpha^\lambda$ is integrable for all $\lambda$, and can be
integrated on the simply connected set $\Sigma_R$ to obtain an element of 
$\MC(\Sigma_R,\hh/\hh^0)|_{-1}^1$ which has the required property.

For the global picture, one can show that, even though a global frame $F$ may not
exist, the equivalence class in $\MC(M,\hh/\hh^0)|_{-1}^1$ is nevertheless well defined.
This argument is given in an analogous situation in \cite{brander2}.
\end{enumerate}
\end{proof}
%****************

%*************************************************************************
%*************************************************************************
\subsection{The Projection to $\Ubar/\Khat$}
We now consider the projection of an element of  $\MC(M,\hh/\hh^0)_{-1}^1$,
to $\Ubar/\Khat$. We assume that $\Ubar$, $\tbar$, $\shat$  are as in Section \ref{constructionsection}
\begin{proposition} \label{flatprop}
Let $F \in \MC(M,\hh/\hh^0)_{-1}^1$.
Let $F_\lambda:M \to \Ubar/(\Kbar\cap \Khat)$ be the associated map, for $\lambda \in \real^*$.
 Then the map $\hat{f}_\lambda : M \to \Ubar/\Khat$
obtained by projection  is a curved flat. 
Moreover, if the projection $f_\lambda: M \to \Ubar/\Kbar$
is regular then so is $\hat{f}_\lambda$.
\end{proposition}
\begin{proof}
Choosing a frame, $F \in \MC(\Sigma,\hh)_{-1}^1$, for $F$, on some contractible
open set $\Sigma$,
$F_\lambda$ is also a frame for $\hat{f}_\lambda$, and, we recall
 the expression for the
Maurer-Cartan form 
$\alpha^\lambda= F_\lambda^{-1}\dd F_\lambda =\alpha_0^{++} + \alpha_1^{+-} (\lambda - \lambda^{-1}) +  \alpha_1^{--}(\lambda + \lambda^{-1})$, where the second two terms 
\beq \label{alpha-}
\alpha_{\phat} := \alpha_1^{+-} (\lambda - \lambda^{-1}) +  \alpha_1^{--}(\lambda + \lambda^{-1}),
\eeq
make up the $\phat$ component of $\alpha$. 
Recall that $\hat{f}_\lambda$ is a curved flat if and only if $\alpha_{\phat} \wedge \alpha_{\phat}
= 0$.  But this follows immediately from the conditions (\ref{fourthcond}) and
(\ref{lastcond}) obtained in Lemma \ref{conditionlemma}.
Finally, if $f_{\lambda}:M \to \Ubar/\Kbar$ is regular, then, as discussed above,
the 1-form $\alpha_1^{--}$ consists of $\textup{Dim}(M)$ linearly independent
1-forms. 
From (\ref{alpha-}), it follows that  $\alpha_{\phat}$ has the same property, and $\hat{f}_\lambda$ is
also regular.
\end{proof}

\begin{corollary} \label{flatpropcor}
Let $F \in \MC(M,\hh/\hh^0)_{-1}^1$, where $\textup{Dim}(M) = n$, and
$\Ubar/\Khat$ is Riemannian.
Suppose that either of the 
associated maps $f_\lambda : M \to \Ubar/\Kbar$ or $\hat{f}_\lambda: M \to \Ubar/\Khat$
is regular for some (and hence all)  values of $\lambda \in \real^*$.
 Then $M$ admits a topological covering by $\real^n$. 
\end{corollary}
\begin{proof}
The map $\hat{f}_\lambda$ is a regular curved flat into a Riemannian symmetric space, 
which is intrinsically flat \cite{feruspedit1996II}.
\end{proof}

Combined with the converse part of Theorem \ref{firstprojectionthm} we obtain
Corollary \ref{flatpropcor2}.

%***************************************************
\section{Regularity of Projections  to $\Ubar/\Kbar$}  \label{regularitysection}
We discussed, in Section \ref{construction}, the local construction, 
via the KDPW method, of
 elements of $\MC(M,\hh/\hh^0)_{-1}^1$ which evaluate to
  immersions into
$\Ubar/(\Kbar\cap \Khat)$, 
provided $\textup{dim}(M) \leq \textup{rank}(\Ubar/\Khat)$.
 Now we look at 
the regularity of projections to $\Ubar/\Kbar$. 

\begin{lemma} \label{projectionlemma}
Suppose that the symmetric space $\Ubar/\Khat$ is Riemannian, and let $V$ be 
any vector subspace of $\phat$ of dimension $n \leq r := \textup{Rank}(\Ubar/\Khat)$.
Let $\pi_V: \phat \to V$ be the orthogonal projection with respect to the
Killing metric on $\phat$.
Then there exists a Cartan subalgebra $\mhat$ of $\phat$ such that 
\bdm
\pi_V \mhat = V.
\edm
Moreover, the set of Cartan subalgebras of $\phat$ which have this property is open and
dense in the Grassmannian of Cartan subalgebras of $\phat$.
\end{lemma}
\begin{proof}
  Let $\mhat$ be any
Cartan subalgebra of $\phat$, and choose an orthonormal basis $\{a_1,...,a_r \}$ consisting 
of regular elements, 
ordered such that the projection 
$\pi_V \text{Span}\{a_1,...,a_j\} = \pi_V \mhat =: V_1$,
for some $j \leq r$.
If the projection is not surjective,
 let $y$ be any non-zero element of $V \cap \mhat^\perp$.  
Then, clearly $y \notin \mhat$, so, since a Cartan subalgebra of $\phat$ is equal
to the centralizer of any of its regular elements, we must have $[a_i,y] \neq 0$
for any $i$.  Then
\bdm
\langle [[a_{j+1},y],a_{j+1}]\, , \, y \rangle = \langle [a_{j+1},y]\, , \, [a_{j+1},y] \rangle \neq 0,
\edm
so the  action of $\phi_t := \Ad_{\exp(t[a_{j+1},y])}$ moves $a_{j+1}$ a non-zero
amount in the direction of $y$, to first order in $t$.  On the other hand, for all $i\neq j+1$,
since $a_k$ are orthonormal, we have
\bdm
0 = \langle a_{j+1} \, , \, a_i\rangle = \langle[a_{j+1},y] \, , \, [a_i,y]\rangle = \langle[[a_{j+1},y],a_i] \, ,\, y\rangle,
\edm
so  $\phi_t a_{i}$ has no first order component in the direction of $y$.   Choose $t$
 small enough
so that $\pi_V  \textup{Span}\{\phi_t a_1,...,\phi_t a_j\}$ still has dimension $j$.
Then
\bdm
\pi_V \phi_t a_{j+1} = Ct y + o(t),
\edm
with $C \neq 0$, and, for $i \neq j+1$,
\bdm
\pi_V \phi_t a_{i} = \pi_V a_i + o(t)w + o(t^2)y,
\edm
with $\langle w \, ,\, y \rangle =0$, and $\langle \pi_V a_i \, , \, y\rangle = 0$.
It follows that $\pi_V \phi_t a_{j+1}$ is linearly independent from $\pi_V \phi_t a_{i}$.
Hence $\pi_V \textup{Span}\{\phi_t a_1,...,\phi_t a_{j+1} \} =: V_2$ has dimension
$j+1$.  

Now  set $\mhat_2 := \phi_t \mhat$. This
 is still a Cartan subalgebra of $\phat$, and $\phi_t a_i$ are still regular
elements, since the adjoint action of $\Khat$ preserves these properties.  Hence, if $j+1<n$,
 we can repeat the above procedure, choosing an element $y_1 \in V \cap \mhat_2^\perp$.
After $n-j$ steps, we have $\pi_V \mhat_{n-j} = V$.  

The above argument showed that there is a Cartan subalgebra which projects surjectively
onto $V$ arbitrarily close to any Cartan subalgebra of $\mhat$ of $\phat$. Since this
projection property is also an open condition among Cartan subalgebras, this means
that a generic Cartan subalgebra projects onto $V$.
\end{proof}

Let $r := \textup{Rank}(\Ubar/\Khat)$.
Regarding the next proposition, we first observe that a regular,
$r$-dimensional, family 
of curved flats in  $\Ubar/\Khat$, given by 
$F_+ \in \MC(\real^r, \Lambda G_{\shat \rho})_1^1$,
 certainly exists. 
One could take $F_+$ to be the family associated to an embedding of a flat totally
geodesic submanifold of $\Ubar/\Khat$, for example.  
Further, if $\Ubar/\Khat$ is Riemannian, then the tangent space at 
each point to an $r$-dimensional regular curved flat is a Cartan subalgebra of $\phat$
\cite{feruspedit1996II}.
%*********************************
\begin{proposition}   \label{immersionprop}
Suppose that the symmetric space $\Ubar/\Khat$ is Riemannian and of rank $r$.
\begin{enumerate}
\item[(i)]
If there exists a local solution $F \in \MC(\Sigma,\hh/\hh^0)_{-1}^1$
such that the associated maps $f_{\lambda}: \Sigma \to \Ubar/\Kbar$, for 
$\lambda \in \real^*$, are regular, then 
$n \leq r$.
\item [(ii)]
Conversely, if $n \leq r$ then such a local regular solution
 can be constructed from any regular curved flat of dimension $r$, and 
 the KDPW correspondence.
\end{enumerate}
\end{proposition}
\begin{proof}
\begin{enumerate}
\item[(i)]
First we examine the relationship between an element $F \in \MC(\Sigma,\hh/\hh^0)_{-1}^1$
and the corresponding family of curved flats, $F_+ \in \MC(\Sigma,\Lambda G_{\shat \rho})_1^1$,
obtained via the KDPW correspondence (\ref{taucor}).  They are related by the equation
\bdm
F = F_+G_-,
\edm
for some $G_- \in \MC(\Sigma,\Lambda G_{\shat \rho})_{-\infty}^0$.
The Maurer-Cartan form of $F_+$ has the expansion
\bdm
F_+^{-1}\dd F_+ = \psi \lambda,
\edm
for some $\psi \in \phat \otimes \Omega(\Sigma)$, and $G_-$ has the expansion
\bdm
G_- = \sum_{-\infty}^0 C_i \lambda^i, \hspace{1cm} C_{2j} \in \Khat, ~ C_{2j+1} \in U_-.
\edm
Now
\beqas
F^{-1} \dd F &=& \alpha_0^{++} + \alpha_1^{+-} (\lambda - \lambda^{-1}) +  \alpha_1^{--}(\lambda + \lambda^{-1}) \\
&=&G_-^{-1}(F_+^{-1} \dd F_+)G_- + G_-^{-1}\dd G_-\\
&=& C_0^{-1} (\psi \lambda) C_0 + \sum_{i=-\infty}^0 D_i \lambda^i,
\eeqas
where $D_i$ and $C_0$ do not depend on $\lambda$.  As discussed above, the coframe for the
projection to $\Ubar/\Kbar$ is obtained from the 1-form $\alpha_1^{--}(\lambda+\lambda^{-1})$, and
so the condition for $f: \Sigma \to \Ubar/\Kbar$ to be an immersion for all $\lambda \in \real^*$
is that $\alpha_1^{--}$ consists of $n$ linearly independent 1-forms.
From the above expansion of $F^{-1}\dd F$, we have 
\beq \label{alpha11eqn}
\alpha_1^{--} = \pi_{\pp} \Ad_{C_0^{-1}} \, \psi.
\eeq
Suppose first that we have a solution such that $f_\lambda$ is an immersion.
Because $F_+$ is a curved flat in $\Ubar/\Khat$, it is necessarily tangent to a flat of $\Ubar/\Khat$,
which means that $\psi$ takes its values in an Abelian subalgebra $\ahat_1 \subset \phat$.
Since the adjoint action of $\Khat$ on $\phat$ takes Abelian subalgebras to Abelian
subalgebras, $\Ad_{C_0^{-1}} \, \psi$ also takes it values in an Abelian subalgebra 
$\ahat \subset \phat$.  If $\pi_{\pp} \Ad_{C_0^{-1}} \, \psi$ consists of $n$ linearly
independent 1-forms, then it is necessary  that $\pi_{\pp} \ahat$ has dimension $n$,
and so the rank of $\Ubar/\Khat$ is at least $n$.\\

\item[(ii)]
%----------------------
Conversely, suppose that $n \leq r$.
By Lemma \ref{projectionlemma}, there exists a Cartan subalgebra
 $\ahat$ of $\phat$ such that $\pi_{\pp}\ahat = \pp$. 
Take any curved flat family $F_+ \in \MC(\Sigma, \Lambda G_{\shat \rho})_1^1$
which is regular, that is,
 $\psi \lambda|_x = F_+^{-1}\dd F_+|_x$ consists of $n$ linearly independent 1-forms  and takes its values in a Cartan subalgebra $\mhat_x \subset \phat$,
 for all $x$ in some neighbourhood of the initial condition $x_0 \in \Sigma$.

We can assume that $x_0$ is the point at which all our loop group maps are normalized,
that is, $F(x_0) = F_+(x_0) = G_-(x_0) = I$ for all $\lambda$, where 
$F \in \MC(\Sigma,\hh/\hh^0)|_{-1}^1$ is obtained via the KDPW correspondence (\ref{taucor}).
 We can also assume, after an action by $\Khat$, that $\mhat_{x_0} =\ahat$, 
 since $\ahat$ is Cartan and $\Khat$ acts transitively on the 
 Cartan subalgebras in $\phat$.
At the normalization point $x_0$,  equation (\ref{alpha11eqn}) reduces to
$\alpha_1^{--}|_{x_0} = \pi_{\pp} \, \psi|_{x_0}$. Since, $\pi_{\pp} \ahat = \pp$, and $\ahat$ and $\pp$ have the same dimension,
 it follows that  $\alpha_1^{--}|_{x_0}$ also consists of $n$ linearly independent
1-forms.  This is an open condition, and hence the solution is an immersion on some neighbourhood
of $x_0$.
\end{enumerate}
\end{proof}

Combining Theorem \ref{firstprojectionthm} with Proposition \ref{immersionprop} we 
obtain Theorem \ref{theorem1}.

%*************************************************************************
%********************************************************************
\section{The Solutions Associated to Reflective Submanifolds of Compact Symmetric Spaces} \label{compactsection}
In this section we determine the answer to Problem \ref{problem} for the case
that $\Ubar/\Kbar$ is a compact, simply connected, irreducible symmetric space.
Clearly, we need only consider strongly curvature invariant $\VV$-geometries
which admit non-totally geodesic submanifolds, so it is enough to go through
Naitoh's list, given in the introduction.
%--------------------------------------
\subsection{The Geometry of $k$-Dimensional Submanifolds of the Sphere $S^n$, $0<k<n$}
The reflective submanifolds of the sphere $S^n$ are just the totally geodesic
submanifolds, namely spheres of dimension $k<n$. For these, the normal bundle is
flat, so isometric normal curvature preserving $\VV_\pp$-immersions of $N_R$ 
will be isometric immersions with flat normal bundle of a sphere of radius
$\sqrt{R}$.

To check whether solutions exist, 
set $\ubar = \mathfrak{so}(n+1)$, 
$\tbar = \textup{diag}(I_n,-1)$,
where $I_j$ denotes the $j\times j$ identity matrix, $\shat = \textup{diag}(I_k,-I_{n+1-k})$.
Then $\kbar = \mathfrak{so}(n)$, $\khat = \mathfrak{so}(k) \times \mathfrak{so}(n+1-k)$, 
and
\bdm
\pp = \pbar \cap \phat = \bbar{ccc} 0 & 0 & *_{k \times 1} \\ 0 & 0 & 0 \\ *_{1 \times k} & 0 & 0 \ebar,
\edm
where $*_{i\times j}$ indicates an $i\times j$ submatrix.
To identify the associated totally geodesic submanifold $\NP$, we can just
look at $\pp$, which must be the tangent space to the origin. We can take
$\kp = [\pp,\pp]$, and we see that $\NP = SO(k+1)/SO(k) = S^k$.
The rank of $\Ubar/\Khat$ is min$(k, n+1-k)$, so, by Theorem \ref{theorem1},
we need $k \leq n+1-k$. In other word, solutions exist if and only if
$k \leq (n+1)/2$.
%------------------------------------
\subsection{The Geometry of $k$-Dimensional Complex Submanifolds of the Complex Projective
Space $\cc P^n$, $0<k<n$}
This is the analogue of the case of the sphere, substituting 
$\mathfrak{su}(n+1)$ for $\mathfrak{so}(n+1)$. In this case, $N = SU(k+1)/S(U(k)\times U(1))$
and $\Ubar/\Khat = SU(n+1)/(S(U(k) \times U(n+1-k))$.  The dimension of $N$ is $2k$,
and the rank of $\Ubar/\Khat$ is min$(k,n+1-k)$ so there is no solution for any
$k$ or $n$.

%-------------------------------------------------------------
\subsection{The Geometry of $n$-Dimensional Totally Real Submanifolds of $\cc P^n$}
This is the case $\ubar = \mathfrak{su}(n+1)$, $\kbar = \mathfrak{s}(\mathfrak{u}(n)\times \mathfrak{u}(1))$,
 $\khat = \mathfrak{so}(n+1)$.  Then $U/K = \real P^n$, and has dimension $n$,
 and $\Ubar/\Khat$ has rank $n$. Therefore, solutions exist for every $n$.
The $\VV_\pp$-immersions are totally real isometric immersions of $S_R^n$ into
$\cc P^n$.   In this case, the normal bundle is isomorphic to the tangent bundle 
 via the complex structure, so the requirement that the immersion be normal
curvature preserving is contained in the totally real condition.

It may be of interest to see the details explicitly:
consider $S^{2n+1}$ as the unit sphere in $\cc ^{n+1}$, with the Sasakian 
structure on $T S^{2n+1}$ induced by multiplication by $i \in \cc^{n+1}$,
which we denote by $J$. We will call a submanifold $M$ of $S^{2n+1}$
\emph{totally real} if it is totally real as a submanifold of $\cc^{n+1}$, 
which means $JT_xM \subset T^\perp_xM \oplus \real x$ for all $x$.

The Hopf fibration, $\pi: S^{2n+1} \to \cc P^n$,
is given by the projection to $S^{2n+1}/U(1) = \cc P^n$, where the action
of $U(1)$ is by following the integral curves of the vector field $Jx$, for
$x \in S^{2n+1}$. The vertical space at $x$ to the fibration is given by
the line through the origin in $T_x S^{2n+1}$ spanned by $Jx$.
If $M= M^n$ is an $n$-dimensional immersed submanifold of $S^{2n+1}$ which is transverse to the Hopf fibration for all $x$, then the projection to $\cc P^n$, $\pi (M)$ is also
an immersed submanifold of the same dimension.  If $M$ is totally real in
$S^{2n+1}$, then
$\pi (M)$ is also totally real (or \emph{Lagrangian}) in  $\cc P^n$, that is
the complex structure on the tangent bundle of $\cc P^n$, which we also
denote by $J$, takes the tangent space $T_{\pi (x)} \pi (M)$ isomorphically to
its orthogonal complement $T_{\pi (x)}^\perp \pi (M)$.  This follows from the
fact that $J$ is induced from the complex structure on $\cc ^{n+1}$.
Conversely, one can show, \cite{reckziegel}, that if $M$ is simply connected
and $\hat{f}: M \to \cc P^n$ is a Lagrangian immersion then there exists a lift
to a Legendrian, (a special case of our definition of totally real) immersion $f: M \to S^{2n+1}$ such that $\pi f = \hat{f}$.

Let $\Ubar = SU(n+1)$, represented by the matrix subgroup of $SO(2n+2)$ consisting
of all matrices of the form $\tiny{\bbar {cc} A & -B\\ B & A \ebar}$, where $A$ and
$B$ are $(n+1) \times (n+1)$, and such that $\det(A+iB) = 1$. 
The Lie algebra $\ub = su(n+1)$ is thus represented as the subalgebra of $so(2n+2)$
consisting of matrices $\tiny{\bbar{cc} a & b \\ -b & a \ebar}$, where $a$ is skew symmetric,
$b$ is symmetric and $\textup{Tr}(b) = 0$.

For a matrix $F \in \Ubar$, we regard the columns
of $F$ as vectors $[a,b]^t = a + bi \in \real^{2n+2} = \cc^{n+1}$. Then 
$F = [X, JX]$, where $J = \tiny{\bbar {cc} 0 & -I_{n+1} \\ I_{n+1} & 0 \ebar}$
is the complex structure on $\cc^{n+1}$.
An \emph{adapted frame} for a totally real immersion $f:M^n \to S^{2n+1}$
is a map $F: M \to U$ given by 
\beq \label{trframe}
F = [X, f, JX, Jf],
\eeq
where the $n$ columns of $X$ are all in the normal space to the immersion.
Thus the tangent space is contained in the span of the vectors making up
$JX$ and $Jf$. 

 Let $\Ubar = SU(n+1)$ as above.  Let $\shat = \Ad _P$, where
 $P ={\textup{diag}(I_{n+1},-I_{n+1})}$, and $\tbar = \Ad_Q$, for 
 $Q = \textup{diag}(I_n, -1, I_n, -1)$.  Then $\Ubar_{\shat}$ and
 $\Ubar_\tbar$ are isomorphic to  $SO(n+1)$ and $S(U_n \times U_1)$ respectively.

Let $G = U^\cc$ and, as before, define $\shat$,
$\tbar$ and $\rho$ by the extensions  (\ref{symshat}) and (\ref{symtau}),
and the extension (\ref{symrho}) for the complex conjugation on $G$,
and set $\hh = \Lambda G_{\rho \shat \tbar}$. Let $F$ be an
element of $\MC(M,\hh/\hh^0)_{-1}^1$.
Now $\hh^0 = \kp =\kb \cap \khat$ consists of matrices in $su(n+1)$ of the form
$\textup{diag}(*_{n\times n}, *_{1 \times 1}, *_{n \times n}, *_{1 \times 1})$,
and so right multiplication by $\hh^0$ fixes the $(n+1)$'th and last columns
of $F$.   We will therefore regard the
$(n+1)$'th column of $F$ as a map $f$ into $S^{2n+1}$.
The Maurer-Cartan form of $F$ has
the expression
\beq \label{slmcf}
F^{-1} \dd F = \bbar {cccc} \omega & 0 & -(\lambda - \lambda^{-1}) \beta
   & -(\lambda + \lambda^{-1}) \theta \\
   0 & 0 & -(\lambda + \lambda)^{-1} \theta^t & -(\lambda-\lambda^{-1})\alpha \\
   (\lambda-\lambda^{-1}) \beta^t & (\lambda + \lambda^{-1})\theta 
     & \omega & 0 \\
     (\lambda+\lambda^{-1})\theta^t & (\lambda-\lambda^{-1})\alpha & 0&0 
   \ebar,
\eeq
where  $\omega$ and $\beta$ are an $n \times n$ matrix-valued 1-forms, $\theta$
is $n \times 1$ and $\alpha$ is $1 \times 1$. 
The condition (\ref{lemcond}) in this case reduces to the equation
\beq \label{cccond}
\dd \omega + \omega \wedge \omega = 4 \theta \wedge \theta^t,
\eeq
which we shall return to below.

 If we write $F = [X, f, JX, Jf]$,
then we see from (\ref{slmcf}) that $X^t \dd f = 0$, and so the columns of 
$X$ are all contained in the normal space to the image of $f$. Thus, $f$ is
totally real. The coframe for $f$ is given by the column vector 
\bdm
[JX,Jf]^t\, \dd f = [(\lambda+\lambda^{-1})\theta^t, \,(\lambda-\lambda^{-1})\alpha]^t,
\edm
and so the condition that $f$ be an immersion is that the 1-forms which make up
the components of this matrix span a space of dimension $n$.

At $\lambda = 1$, we have $(Jf)^t \, \dd f = (\lambda -\lambda^{-1}) \alpha = 0$,
 so $Jf$ is normal, and the immersion
is Legendrian.  At this point, the second fundamental form is 
\bdm
[X,  Jf]^t \, \dd (JX) = [-(\lambda-\lambda^{-1})\beta, 0]^t,
\edm
and this is also zero, so the immersion is  totally geodesic Legendrian at $\lambda =1$.
This projects to a totally geodesic Lagrangian immersion in $\cc P^n$.

For other values of $\lambda$, $f$ is not horizontal, but, if we assume that
the vector $\theta^t$ consists of $n$ linearly independent 1-forms, then $f$ is
transverse to the fibre of the Hopf projection.  Thus $\pi \circ f$ is still
a Lagrangian immersion into $\cc P^n$.  Writing $X = [X_1,...,X_n]$, 
then, since $X_i$ are horizontal to $\pi$, and normal to the image of $f$,
it follows that $\pi_* X_i$ are an orthonormal basis for the normal space to $\pi f$,
and, therefore, $\pi_* JX_i = J \pi_* X_i$ are an orthonormal basis for the tangent space.

Let $e_i$ be the vector fields on $M$ such that $\pi_* JX_i = \pi_* \circ f_* e_i$. Then
$f_* e_i = JX_i + (f_*e_i)^\perp$, where the last term is in the vertical space.
Hence $JX_i^t \, \dd f (e_j) = \delta_{ij}$, and so the coframe for the orthonormal
 (with respect to the metric induced by $\pi \circ f$) basis $e_i$ is given by
\bdm
\theta_{\lambda} := (JX)^t \dd f = (\lambda + \lambda^{-1}) \theta.
\edm
The connection induced by $\pi \circ f$ on the tangent space is given by
$\omega = (JX)^t \, \dd (JX)$, and  thus the equation (\ref{cccond}), which
we can write as
\bdm
\dd \omega + \omega \wedge \omega = \frac{4}{(\lambda+\lambda^{-1})^2} \theta_{\lambda} \wedge \theta_{\lambda}^t,
\edm
is the statement that the induced curvature is constant and equal to
$c_\lambda = \frac{4}{(\lambda+\lambda^{-1})^2}$. Equivalently, 
 $\pi \circ f$ is an isometric immersion of a piece of 
a sphere of radius $\sqrt{R_\lambda}$.

%----------------------------------------------------------------------------------

\subsection{The Geometry of $2n$-Real-Dimensional Totally Complex Submanifolds of the Quaternionic Projective Space $\HH P^n$.}
The case $n=1$ is a 4-sphere of constant curvature $4$, covered above, so we 
consider $n \geq 2$.
In this case, $\ubar = \mathfrak{sp}(n+1)$, $\kbar = \mathfrak{sp}(n) \times \mathfrak{sp}(1)$,
$\khat = \mathfrak{u}(n+1)$, and the reflective submanifold is $\cc P^n$, which has
dimension $2n$. The rank of $\Ubar/\Khat$ is $n+1$, so we need $2n \leq n+1$, and
solutions do not exist for $n >1$.  

\subsection{The Geometries Associated to Symmetric $R$-Spaces} \label{rspacessubsect}
Naitoh proved that the
  strongly curvature $\VV$-geometries in compact symmetric spaces
  associated to symmetric $R$-spaces  only
  admit non-totally geodesic $\VV$-submanifolds which are homothetic to
 the corresponding reflective submanifold, with a factor $R<1$, opposite to
 what we seek.
  The analogue for the non-compact
case, with $R>1$, was proved by Berndt et al. in \cite{bent} (where both cases are discussed).
 Hence we do not need to consider
these cases.  

Taking into account Corollary \ref{flatpropcor2}, this completes the 
proof of Theorem \ref{compactthm}.

%***********************************************************************

\section{The Non-Compact Case}  \label{noncompactsection}
Here we want to consider the case $R<1$.  We start with the observation 
that the arguments given on the interpretation of the solutions of
$\MC(M,\hh/\hh^0)^1_{-1}$ as projections to $\Ubar/\Kbar$ 
did not depend on what kind of conjugation 
$\rho$ we used for the reality condition.
  The reality condition only determined which subset of $\cc^*$
corresponded to solutions in $\Ubar$, and hence the range of  values for
 $R_\lambda = \frac{(\lambda+\lambda^{-1})^2}{4}$.  

We therefore consider now an alternative reality condition defined by the
conjugation
\bdm
(\rho_2 X)(\lambda) := \rho (X(1/\bar{\lambda})).
\edm
Elements of $\Lambda G_{\rho_2}$ are $\Ubar$-valued for values of $\lambda$ in $S^1$,
and it is easy to check the following
%----------------------
\begin{lemma} \label{s1projectionlemma}
Let $\hh_2 := \Lambda G_{\tbar,\shat,\rho_2}$, where $\shat$ and $\tbar$
are as before. Then  Theorem \ref{firstprojectionthm} also holds
for $\MC(M,\frac{\hh_2}{\hh_2^0})_{-1}^1$, replacing the statement $\lambda \in \real_+$
with $\lambda \in S^1 \setminus \{ \pm i \}$.
\end{lemma}
The values $\lambda = \pm i$ are excluded because the coframe $\theta^\lambda$ vanishes 
there. Note that the difference here is that for $\lambda \in S^1$, $R_{\lambda}$ takes its values in $(0,1)$, and this is the case we are interested in here.

Now define another reality condition $\rt$ by the involution:
\bdm
(\rt X)(\lambda) := \rho \circ \tbar \circ \shat (X(\bar{\lambda})).
\edm
Then $X$ is fixed by $\tbar$, $\shat$ and $\rho_2$ if and only if $X$ 
is fixed by
$\tbar$, $\shat$ and $\rt$.  Hence 
\bdm
\Lambda G_{\tbar \shat \rho_2} = \Lambda G_{\tbar \shat \rt},
\edm
and thus the  loop group is exactly the type we have already considered,
since $\rt$ is of the form (\ref{symrho}). Thus we have the following
%-------------------------
\begin{lemma} \label{rprojectionlemma}
Let $\hh_2 := \Lambda G_{\tbar \shat \rho_2}$, as in Lemma \ref{s1projectionlemma}.
 Then  Theorem \ref{firstprojectionthm} also holds
for $\MC(M,\frac{\hh_2}{\hh_2^0})_{-1}^1$, for $\lambda \in \real_+$,
replacing the pairwise symmetric Lie algebra $(\ubar, \tbar, \shat)$ with
$(\utilde, \tbar, \shat)$,
where $\utilde$ is the real form of the Lie algebra $\mathfrak{g}$ defined by
\bdm
\rt := \rho \circ \tbar \circ \shat.
\edm
\end{lemma}
Hence the case $R<1$ is really just the case $R>1$ for a symmetric space associated
to a different real form of $G$, and so we can use our previous analysis to discuss
the problem of existence of solutions.
Specifically, the projection $f_\lambda : M \to \Ubar/\Kbar$ has solutions
for $\lambda \in S^1$ if and only if the projection $\tilde{f}: M \to \Utilde/\Ktilde$
has solutions for $\lambda \in \real^*$. Set
 $\Kcheck := \Utilde_{\shat}$.
The point that we will use is: if $\Utilde/\Kcheck$ is
 \emph{Riemannian} then solutions exist if and only if
  $\textup{Dim} (U/K) = \textup{Rank}(\Utilde/\Kcheck)$.

%--------------------------------------------
\subsection{The Geometry of $k$-Dimensional Submanifolds of the Real Hyperbolic Space}
This is the analogue of the case of the sphere. Normal curvature preserving $\VV_\pp$-immersions of $N_R$ 
will be isometric immersions with flat normal bundle of a hyperbolic space of radius
$\sqrt{R}$.
We have the analogous construction with $\ubar = \mathfrak{so}(n,1)$,
$\tbar = \Ad_Q$, $Q = \textup{diag}(I_n,-1)$,
 $\shat = \Ad_P$, $P = \textup{diag}(I_k,-I_{n+1-k})$.

Set $J_{p,q} = \textup{diag}(I_p,-I_q)$, and define 
\beqas
SO(p,q,\cc) = \{X \in GL(p+q,\cc) ~|~ X^tJ_{p,q}X = J \},
     \\
SO(p,q) = \{ X \in SO(p,q,\cc)  ~|~ \overline{X} = X \}.
\eeqas
 According to the argument presented
above, if we replace the involution $\rho X = \overline{X}$ with 
\bdm
\rt X = \overline{\tbar \shat X},
\edm
then we can analyze the corresponding objects in the real form
given by $\Utilde = SO(n,1,\cc)_{\rt}$ for existence
of solutions.  Let $\phi: GL(n+1,\cc) \to GL(n+1,\cc)$ be defined by $\phi = Ad_T$,
for $T = \textup{diag}(I_k, iI_n,1)$.  Then it is easy to check that
\bdm
\phi: SO(n,1,\cc)_{\rt} \to SO(k,n+1-k)
\edm
 is an isomorphism. $\phi$ commutes with $\shat$, so we have 
\bdm
\utilde = \mathfrak{so}(k,n+1-k), \hspace{1cm} \kcheck = \mathfrak{so}(k) \times \mathfrak{so}(n+1-k),
\edm
 which is compact. Hence $\Utilde/\Kcheck$ is Riemannian,
 and of rank min$(k,n+1-k)$, 
 and we can say that solutions exist if and only if $k \leq (n+1)/2$.

\begin{remark} The solutions $\tilde{f}_\lambda: M \to \Utilde/\Ktilde$ here, for
$\lambda \in \real^*$, are actually strongly curvature invariant
$\VV$-immersions, $H^k_R \to SO(k,n+1-k)/SO(k,n-k)$, and do not appear in
Naitoh's list because the target space is non-Riemannian. 
\end{remark} 
%---------------------------------------------------
\subsection{The Geometry of Complex Submanifolds of Complex Hyperbolic Space $\cc H^n$}
This is the  analogue of the previous case, substituting 
$\mathfrak{su}(n,1)$ for $\mathfrak{so}(n,1)$. Then
$\utilde = \mathfrak{su}(k,n+1-k)$, and
$\kcheck = \mathfrak{s}(\mathfrak{u}(k) \times \mathfrak{u}(n+1-k))$, 
so $\Utilde/\Kcheck$ is Riemannian, and, as in the compact analogue,
there is no solution for any $k$ or $n$.
%------------------------------------------------
\subsection{The Geometry of $n$-Dimensional Totally Real Submanifolds of $\cc H^n$}
Here $\ubar = \mathfrak{su}(n,1)$, $\kbar = \mathfrak{s}(\mathfrak{u}(n)\times \mathfrak{u}(1))$,
 $\khat = \mathfrak{so}(n,1) $. The complex Lie algebra is $\mathfrak{g} = \mathfrak{sl}(n+1,\cc)$, and the involutions on $\mathfrak{g}$ are given by 
 $\rho x = \Ad_J (\bar{x}^t)^{-1}$, $\tbar x = \Ad_J x$ and $\shat x = \Ad_J (x^t)^{-1}$,
 where $J = \textup{diag}(I_n, -1)$.
 
 Setting $\rt (x) = \rho \circ \tbar \circ \shat (x) = \Ad_J \bar{x}$,
 we obtain
\bdm
\tilde{\mathfrak{u}} = \mathfrak{sl}(n+1,\real), \hspace{1cm} 
\kcheck = \mathfrak{so}(n,1).
\edm
In this case, $\kcheck$ is not compact, so $\widetilde{U}/\Kcheck$ is non-Riemannian.
However, we do not want to show non-existence
here, so this is no problem.  To show that a local solution \emph{can} be constructed from a
curved flat in $\widetilde{U}/\Kcheck$, via the KDPW scheme,
examining the converse part of the proof of Proposition \ref{immersionprop},
it will be enough to find an Abelian subalgebra $\mathfrak{a}$ 
in $\utilde_-$ which projects
onto $\ptilde \cap \utilde_-$.  Then exponentiating $\mathfrak{a}$ 
will give a suitable curved flat to locally construct a regular solution, 
as in the proof of Proposition \ref{immersionprop}.
  The relevant subspaces consist of the matrices in $\mathfrak{sl}(n+1,\real)$
as follows:
\bdm
\utilde_- = \left\{ \bbar {cc} A & b\\ -b^t & c \ebar \right\}, \hspace{1cm}
\ptilde \cap \utilde_-   \left\{ \bbar {cc}  0 & *_{n \times 1} \\  
      *_{1 \times n} & 0 \ebar \right\},
\edm
where $A$ is a symmetric $n \times n$ matrix.  Let $E_i$ be the $(n+1) \times (n+1)$
matrix  whose $i$'th row is $[0,...,0,1,1]$ and $i$'th column is $[0,...,0,1,-1]$,
and all other components are zero. Let $\hat{E}$ be the matrix whose $n$'th row is
 $[0,...,0,1,1]$, and whose $(n+1)$'st row is $[0,...,0,-1,-1]$, with all other
 components zero.
 Then $\{E_1,...,E_{n-1},\hat{E}\}$ generate such an Abelian subalgebra.
Hence, a solution exists for every $n$.

\subsection{The Geometry of $2n$-Real-Dimensional Totally Complex Submanifolds of the Quaternionic Hyperbolic Space $\HH H^n$}

In this case, $\mathfrak{g} = \mathfrak{sp}(n+1,\cc)$, 
$\ubar = \mathfrak{sp}(n,1)$, $\kbar = \mathfrak{sp}(n) \times \mathfrak{sp}(1)$,
$\khat = \mathfrak{u}(n,1)$, and the reflective submanifold is $\cc H^n$, which has
dimension $2n$. 

The involutions are $\rho X = \Ad_K \overline{(X^{t})^{-1}}$, 
for $K = \textup{diag}(I_n,-1,I_n,-1)$, $\tbar X = \Ad_K X$, $\shat X = \Ad_J X$, for
$J = \tiny{\bbar {cc} 0 & I_{n+1} \\ -I_{n+1} & 0 \ebar}$.
Setting $\rt (X) = \rho \circ \tbar \circ \shat (X) =  \Ad_J \overline{(X^t)^{-1}} = \bar{X}$,
we have
\bdm
\tilde{\mathfrak{u}} = \mathfrak{sp}(n+1,\real), \hspace{1cm} 
  \check{\mathfrak{k}} = \mathfrak{u}(n+1).
\edm
Now $\widetilde{U}/\Kcheck$ is  Riemannian and of rank $n+1$, hence,
 as in the compact case,
 we need $2n \leq n+1$, so 
solutions do not exist for  $n>1$.  

Together with the remarks of Section \ref{rspacessubsect}, this completes
the proof of Theorem \ref{noncompactthm}.

%-----------------------------------------------------

%***********************************************
\section{Projections to Other Homogeneous Spaces}
     \label{surfacedescription}

There are many more reflective submanifolds besides those discussed in 
Sections \ref{compactsection} and \ref{noncompactsection}.  In all other 
(Riemannian) cases, the rank $r$ of $\Ubar/\Khat$ must be less than the dimension 
of $\pp$, according to the results of Naitoh.
  However, we can still produce solutions of dimension $r$, and 
 these correspond, under the projection to $\Ubar/\Kbar$ to
 \emph{sub}-$\VV_\pp$-immersions. In this case,
they have the property that the connections on the natural bundles $\pp$ and 
$\pp^\perp$ are preserved, as $\lambda$ varies,
 rather than the tangent and normal bundles.

It may be be interesting to take projections to other homogeneous spaces
$\Ubar/H$ where $H$ is a subgroup containing $\Kbar \cap \Khat$, and
we consider one example here.

\subsection{Reflective Submanifolds of Symmetric Spaces with Isometry Group $G_2$}
According to the classification in \cite{leung1974}, the only
reflective submanifolds associated to $G_2$ are 
the reflective submanifolds $SO(4)$ and $G_2/SO(4)$ in $G_2$, and 
 $(SU(2)/SO(2)) \times (SU(2)/SO(2))$ in $G_2/SO(4)$. All of these
 have dimension greater than $2$, but the rank of $G_2$ is $2$,
and so the rank of $G_2/(G_2)_+$ will always be too small
for the existence of local normal curvature preserving $\VV_\pp$-immersions
of $N_R$, for $R>1$, for any of these reflective submanifolds.

\subsection{Surfaces in $S^6$ with $G_2$-Frames}
Let $i,~ j,~ k,~x,~ ix,~ jx$ and  $kx$ be a basis for the imaginary
octonions, $\imo$, and identify these with $e_1,...,e_7$, the standard basis
for ${\real}^7$.  We define the multiplication on the octonions, 
by the usual quaternionic multiplication on $i,j$ and $k$,
and extending the following table by linearity and the relation $ab=-ba$ for
all imaginary octonions:
\begin{table}[here]  
  \begin{tabular}{|c|c|c|c|c|c|c|c|}  \hline    
   & i & j   & k & x & ix & jx & kx \\ 
 \hline
ix & x & -kx & jx & -i&-1  &-k & j \\
\hline
jx & kx&   x & -ix&-j& k  &-1& -i\\
\hline
kx & -jx& ix & x  & -k &-j & i& -1\\
  \hline 
  \end{tabular} 
\end{table}

Recall that $G_2$ is the Lie group consisting of algebra automorphisms 
of $\OO$. Since an element of $G_2$ fixes the identity, it is determined
by its (orthogonal, orientation preserving) action on $\imo = \real ^7$;
therefore, $G_2$ is a subgroup of $SO(7)$. A \emph{$G_2$-frame} is
a matrix $F = [f_1,...,f_7] \subset G_2 \subset SO(7)$, whose columns
$f_1,...,f_7,$ necessarily satisfy the the same octonionic 
multiplication table as $e_1$,...,$e_7$, described above. That is,
$f_1 f_2 = f_3$, and so on.

There is an almost complex structure $J$ on $S^6 \subset
{\real}^7 = \imo$, defined by Calabi, 
\cite{calabi1958}, as follows: if $u \in S^6$, one can identify the 
tangent space $T_uS^6$ with the 6-plane in ${\real}^7$ orthogonal to
$u$. Then $J_u:~ T_uS^6 \to T_uS^6$ is defined by right multiplication
by $u$. This is an (orthogonal) linear transformation which satisfies 
$J_u ^2 = -1$, and therefore an almost complex structure. 

\subsection{Complex Curves in $S^6$}
A surface in $S^6$ whose tangent space is invariant under the almost complex
structure described above is called an \emph{(almost) complex curve} in $S^6$. 
 If $f: M \to S^6$ is an immersed complex curve,
one way to choose an adapted $G_2$-frame, $F: M \to G_2$ is,
\beq \label{acframe}
F = [f, n_1, n_2, n_3, n_4, e_1, e_2],
\eeq
where $n_i$ are normal, and $e_i$ are tangent to the surface.
For such a matrix in $G_2$, $e_i f = \pm e_j$.  

A  study of complex curves in $S^6$ can be found in \cite{bryant1982}.
They have been studied recently as an integrable system
associated to primitive maps into a 6-symmetric space $G_2/T^2$, 
in \cite{terngkongwang}.

\subsection{Special Surfaces from the Three Involution Loop Group}

Let us now apply the construction of Section \ref{constructionsection} to 
the group $G_2$, that is, we wish to interpret elements of $\MC(M,\hh/\hh^0)_{-1}^1$,
where $\hh = \Lambda (G_2^\cc)_{\rho \shat \tbar}$, and $\rho$, $\shat$ and
$\tbar$ are defined, respectively,
 by the extensions (\ref{symrho}), of complex conjugation,
and the extensions  (\ref{symshat}) and (\ref{symtau}) of
\beqas
\shat = Ad_P, \hspace{1cm} P = \textup{diag}(I_3, -I_4),\\
\tbar = Ad_Q, \hspace{1cm} Q = \textup{diag}(1, -I_2, I_2, -I_2).
\eeqas

%Here $\kb = (\mathfrak{g}_2)_\tbar$ and $\khat = (\mathfrak{g}_2)_{\shat}$
 % are both isomorphic to $\mathfrak{so}(4)$, $\kp \cong \mathfrak{so}(2) \times \mathfrak{so}(2)$
%and $\kp + \pp \cong \mathfrak{su}(2) \times  \mathfrak{su}(2)$.

The subspaces $\kp = \kb \cap \khat$ and $\pp = \pbar \cap \phat$
consist of the elements of $\mathfrak{g}_2$ which have non-zero entries as 
follows:
\bdm
\kp = \bbar {cccc} 0 & & &\\ & *_{2\times2} & &\\ & & *_{2\times2} & \\ &&&*_{2\times 2}\ebar,
\hspace {1cm} 
\pp =\bbar {cccc} &&&*_{1\times 2} \\ && *_{2\times2}&\\ & *_{2\times2} &&\\ *_{2\times 1} &&&\ebar,
\edm

Since right multiplication by an element of $\KP$ fixes the first column of
a matrix $F \in G_2$, there is a natural projection of an element
$F \in G_2/\KP$ onto $S^6 = G_2/SU(3)$, given by the taking this
column.  Because $\hh^0 = \KP$, 
 an element $F \in \MC(M,\hh/\hh^0)_{-1}^1$ can, for real values of
 $\lambda$, be interpreted,  via this projection, 
 as a $G_2$-frame for a certain map $f:M \to S^6$.
Write the frame as
\bdm
F= [f,N,X,Y] = [f,N_1,N_2,X_1,X_2,Y_1,Y_2],
\edm
where $f$ is a single column and $N$, $X$ and $Y$ consist of two columns each.
The last $6$ columns of $F$ form an orthonormal basis for the total bundle
$TS^6|_M$, and the action by $K$ on the right of such a frame preserves
the 2-dimensional subspaces spanned by $N$, $X$ and $Y$ respectively.
Thus, given that the frame is normalized to the identity at some point,
 there is also a natural decomposition of $TS^6|_M$ into these
three 2-dimensional sub-bundles,
\bdm
TS^6|_M = \eta_1 \oplus \eta_2 \oplus \eta_3.
\edm
 Further, each 
sub-bundle, $\eta_i$, has an almost complex structure on it, inherited from 
$S^6$, because the $\{f,N_1,N_2\}$, $\{f,X_1,X_2\}$ and $\{f,Y_1,Y_2\}$
are all associative triples.

Now if $F \in \MC(M,\hh)_{-1}^1$, then the Maurer-Cartan form of $F$ is of the form
\beq \label{scmcf}
 F^{-1}\dd F = \bbar {cccc} 0 &  0 & -(\lambda - \lambda^{-1})\theta_1^t & 
      -(\lambda + \lambda^{-1})\theta_2^t \\
      0 & \omega_1 &  (\lambda + \lambda^{-1})\beta_1 & 
         (\lambda - \lambda^{-1})\beta_2 \\
         (\lambda - \lambda^{-1})\theta_1 & - (\lambda + \lambda^{-1})\beta_1^t
          & \omega_2 & 0 \\
          (\lambda + \lambda^{-1})\theta_2 & -(\lambda - \lambda^{-1})\beta_2^t
           & 0 & \omega_3 \ebar,
\eeq
where $\beta_i$ and $\omega_i$ are all $2\times 2$ submatrices
and $\theta_i$ are $2\times 1$ column submatrices. 
Noting that 
\beq \label{second}
F^{-1} \dd F = \bbar {cccc} f^t \dd f & f^t \dd N & f^t \dd X & f^t \dd Y \\
        N^t \dd f & N^t \dd N & N^t \dd X & N^t \dd Y \\
        X^t \dd f & X^t \dd N & X^t \dd X & X^t \dd Y \\
        Y^t \dd f & Y^t \dd N & Y^t \dd X & Y^t \dd Y \ebar,
\eeq
it follows from (\ref{scmcf}) that the bundle $\eta_1$ is a sub-bundle of the normal
bundle, that $\omega_1$, $\omega_2$ and $\omega_3$ are the
connections for the bundles $\eta_1$, $\eta_2$ and $\eta_3$ respectively,
and that $[(\lambda - \lambda^{-1})\theta_1, (\lambda + \lambda^{-1})\theta_2]^t$
is essentially the coframe.

Note that $\phat$ consists of all matrices in $\mathfrak{g}_2$ of the form
\bdm
\bbar {cc} 0 & *_{3\times 4} \\ *_{3\times4} & 0 \ebar,
\edm
so it follows from Lemma \ref{projectionlemma} that solutions constructed
from 2-dimensional regular curved flats via the KDPW method will generically 
project to regular maps into $S^6$.

At $\lambda =1$ we have
\bdm
 F^{-1}\dd F \big |_{\lambda =1} = \bbar {cccc} 0 &  0 & 0 & 
      -2 \theta_2^t \\
      0 & \omega_1 &  2 \beta_1 & 
         0 \\
         0 & - 2 \beta_1^t
          & \omega_2 & 0 \\
          2\theta_2 & 0
           & 0 & \omega_3 \ebar.
\edm
Comparing with (\ref{second}), this implies that, for this value of $\lambda$,
$X_1$ and $X_2$ also lie in the normal bundle, so the tangent space is spanned by
$Y_1$ and $Y_2$. Hence, $f|_{\lambda =1}$ is a
 complex curve. Moreover, the second fundamental form for $f|_{\lambda =1}$
is given by the matrix
$[N^t \dd Y, X^t \dd Y]^t$ and this is zero. Thus $f|_{\lambda =1}$ is a totally 
geodesic  complex curve in $S^6$.

We now examine the conditions we get for $f_\lambda$  from the equation
(\ref{lemcond}).  It is easy to verify that, in this case, the condition is
equivalent to the three equations
\beqa
\dd \omega_1 + \omega_1 \wedge \omega_1 = 4 \beta_1 \wedge \beta_1^t,\\
\dd \omega_2 + \omega_2 \wedge \omega_2 = 4 \beta_1^t \wedge \beta_1,\\
\dd \omega_3 + \omega_3 \wedge \omega_3 = 4 \theta_2 \wedge \theta_2^t.
\eeqa
These equations say that the
bundle-connection pairs $(\eta_i, \nabla_i)$, where the connection is induced
from $S^6$, are isomorphic to those obtained at $\lambda = 1$.

\section*{Funding}
Research supported by the Japan Society for the Promotion of Science.

\section*{Acknowledgements}
Thanks are due to Josef Dorfmeister, Daniel Fox,  Jost-Hinrich
Eschenburg and Martin Guest, for helpful advice at various 
critical moments  during the preparation of this manuscript.


\begin{thebibliography}{10}

\bibitem{akns1}
M.~J. Ablowitz, D.~J. Kaup, A.~C. Newell, and H.~Segur.
\newblock Method for solving the sine-{G}ordon equation.
\newblock {\em Phys. Rev. Lett.}, 30:1262--1264, 1973.

\bibitem{adlervanmoerbeke1}
M.~Adler and P.~van Moerbeke.
\newblock Completely integrable systems, {Euclidean} {Lie} algebras, and
  curves.
\newblock {\em Adv. in Math.}, 38:267--317, 1980.

\bibitem{adlervanmoerbeke2}
M.~Adler and P.~van Moerbeke.
\newblock Linearization of {Hamiltonian} systems, {Jacobi} varieties and
  representation theory.
\newblock {\em Adv. in Math.}, 38:318--379, 1980.

\bibitem{bent}
J.~Berndt, J.~H. Eschenburg, H.~Naitoh, and K.~Tsukada.
\newblock Symmetric submanifolds associated with irreducible symmetric
  {R}-spaces.
\newblock {\em Math. Ann.}, 332:721--737, 2005.

\bibitem{brander2}
D.~Brander.
\newblock Curved flats, pluriharmonic maps and constant curvature immersions
  into pseudo-{R}iemannian space forms.
\newblock {\em arxiv preprint:0610058. To appear in: \emph{Ann. Global Anal.
  Geom.}}

\bibitem{branderdorf}
D.~Brander and J.~Dorfmeister.
\newblock The generalized {DPW} method and an application to isometric
  immersions of space forms.
\newblock {\em arxiv preprint:0604247}, 2006.

\bibitem{bdpt2002}
M.~Br\"{u}ck, X.~Du, J.~Park, and C.~L. Terng.
\newblock The submanifold geometries associated to grassmannian systems.
\newblock {\em Mem. Amer. Math. Soc. 155}, (735):viii + 95 pp., 2002.

\bibitem{bryant1982}
R.~L. Bryant.
\newblock Submanifolds and special structures on the octonians.
\newblock {\em J. Differential Geom.}, 17(2):185--232, 1982.

\bibitem{burstallpedit}
F.~E. Burstall and F.~Pedit.
\newblock Harmonic maps via {Adler}-{Kostant}-{Symes} theory.
\newblock In {\em Harmonic maps and integrable systems}, number E23 in Aspects
  of Mathematics. Vieweg, 1994.

\bibitem{calabi1958}
E.~Calabi.
\newblock Construction and properties of some $6$-dimensional almost complex
  manifolds.
\newblock {\em Trans. Amer. Math. Soc.}, 87:407--438, 1958.

\bibitem{dorfmeisterpeditwu}
J.~Dorfmeister, F.~Pedit, and H.~Wu.
\newblock Weierstrass type representation of harmonic maps into symmetric
  spaces.
\newblock {\em Comm. Anal. Geom.}, 6:633--668, 1998.

\bibitem{feruspedit1996II}
D.~Ferus and F.~Pedit.
\newblock Curved flats in symmetric spaces.
\newblock {\em Manuscripta Math.}, 91:445--454, 1996.

\bibitem{feruspedit1996}
D.~Ferus and F.~Pedit.
\newblock Isometric immersions of space forms and soliton theory.
\newblock {\em Math. Ann.}, 305:329--342, 1996.

\bibitem{harveylawson}
R.~Harvey and H.~B. Lawson.
\newblock Calibrated geometries.
\newblock {\em Acta Math.}, 148:47--157, 1982.

\bibitem{hilbert1}
D.~Hilbert.
\newblock Ueber {Flachen} von constanter {Gausscher} {Krummung}.
\newblock {\em Trans. Amer. Math. Soc.}, 2:87--99, 1901.

\bibitem{terngkongwang}
S.~Kong, C.~L. Terng, and E.~Wang.
\newblock Associative cones and integrable systems.
\newblock {\em Chinese Ann. Math. Ser. B}, 27(2):153--168, 2006.

\bibitem{kostant1979}
B.~Kostant.
\newblock The solution to a generalized {Toda} lattice and representation
  theory.
\newblock {\em Adv. in Math.}, 34:195--338, 1979.

\bibitem{krichever}
I.~M. Krichever.
\newblock An analogue of the d'{A}lembert formula for the equations of a
  principal chiral field and the sine-{G}ordon equation.
\newblock {\em Dokl. Akad. Nauk SSSR}, 253(2):288--292, 1980.

\bibitem{leung1974}
D.~S.~P. Leung.
\newblock On the classification of reflective submanifolds of {R}iemannian
  symmetric spaces.
\newblock {\em Indian Univ. Math. J.}, 24:327--339, 1974.

\bibitem{leung1975}
D.~S.~P. Leung.
\newblock Errata: ``{O}n the classification of reflective submanifolds of
  {R}iemannian symmetric spaces''.
\newblock {\em Indian Univ. Math. J.}, 24:1199, 1975.

\bibitem{leung1979}
D.~S.~P. Leung.
\newblock Reflective submanifolds. {III}. congruency of isometric reflective
  submanifolds and corrigenda to the classification of reflective submanifolds.
\newblock {\em J. Differential Geom.}, 14, 1979.

\bibitem{moore1}
J.~D. Moore.
\newblock Isometric immersions of space forms in space forms.
\newblock {\em Pacific J. Math.}, 40:157--166, 1972.

\bibitem{naitoh1986}
H.~Naitoh.
\newblock Symmetric submanifolds of compact symmetric spaces.
\newblock {\em Tsukuba J. Math.}, 10:215--242, 1986.

\bibitem{naitoh1993}
H.~Naitoh.
\newblock Compact simple lie algebras with two involutions and submanifolds of
  compact symmetric spaces {I}, {II}.
\newblock {\em Osaka J. Math.}, 30, 1993.

\bibitem{naitoh1998}
H.~Naitoh.
\newblock Grassman geometries on compact symmetric spaces of general type.
\newblock {\em J. Math. Soc. Japan}, 50, 1998.

\bibitem{naitoh2000II}
H.~Naitoh.
\newblock Grassman geometries on compact symmetric spaces of classical type.
\newblock {\em Japanese J. Math.}, 26, 2000.

\bibitem{naitoh2000I}
H.~Naitoh.
\newblock Grassman geometries on compact symmetric spaces of exceptional type.
\newblock {\em Japanese J. Math.}, 26, 2000.

\bibitem{pressleysegal}
A.~Pressley and G.~Segal.
\newblock {\em Loop Groups}.
\newblock Oxford Math. monographs. Clarendon Press, Oxford, 1986.

\bibitem{reckziegel}
H.~Reckziegel.
\newblock Horizontal lifts of isometric immersions into the bundle space of a
  pseudo-{R}iemannian submersion.
\newblock {\em Global differential geometry and global analysis (1984), Lecture
  Notes Math.}, 1156:264--279, 1985.

\bibitem{symes1980}
W.~W. Symes.
\newblock Systems of {Toda} type, inverse spectral problems, and representation
  theory.
\newblock {\em Invent. Math.}, 59:13--51, 1980.

\bibitem{terngtenenblat}
K.~Tenenblat and C.~L. Terng.
\newblock Baecklund's theorem for $n$-dimensional submanifolds of
  $\real^{2n-1}$.
\newblock {\em Ann. Math.}, 111:477--490, 1980.

\bibitem{terng1980}
C.~L. Terng.
\newblock A higher dimensional generalisation of the sine-{G}ordon equation and
  its soliton theory.
\newblock {\em Ann. Math.}, 111:491--510, 1980.

\bibitem{terng2002}
C.~L. Terng.
\newblock Geometries and symmetries of soliton equations and integrable
  elliptic equations.
\newblock {\em arXiv preprint:0212372}, 2002.

\end{thebibliography}
\end{document}